\documentclass[12pt]{amsart}
\usepackage{amsmath}
\usepackage{amsthm}
\usepackage{amssymb}
\usepackage{mathrsfs}
\usepackage{cite}
\usepackage{amscd}
\usepackage{fancyheadings}
\usepackage[all]{xy}

\theoremstyle{definition}
\newtheorem{thm}{Theorem}[section]
\newtheorem{definition}[thm]{Definition}
\newtheorem{prop}[thm]{Proposition}
\newtheorem{lem}[thm]{Lemma}
\newtheorem{rem}[thm]{Remark}
\newtheorem{cor}[thm]{Corollary}
\newtheorem{ex}[thm]{Example}

\newtheorem*{ack}{Acknowledgement}
\newtheorem{para}{}[subsection]
\numberwithin{equation}{section}

\usepackage[top=20truemm, bottom=20truemm, left=20truemm, right=20truemm]{geometry}

\newcommand{\Z}{\mathbb{Z}}
\newcommand{\Q}{\mathbb{Q}}
\newcommand{\Spec}{\operatorname{Spec}}

\newcommand{\Frac}{\operatorname{Frac}}

\newcommand{\Hom}{{\rm Hom}}
\newcommand{\Coh}{\mathsf{Coh}}
\newcommand{\Rep}{\mathsf{Rep}}
\newcommand{\Vecf}{\mathsf{Vec}}
\newcommand{\Vect}{\mathsf{Vect}}
\newcommand{\Mod}{\mathsf{Mod}}
\newcommand{\Qcoh}{\mathsf{QCoh}}
\newcommand{\EFin}{\mathsf{EFin}}
\newcommand{\TFin}{\mathsf{TFin}}
\newcommand{\Aut}{{\rm Aut}}
\newcommand{\Dim}{{\rm dim}}

\newcommand{\GL}{{\rm GL}}

\newcommand{\Ker}{\operatorname{Ker}}
\newcommand{\im}{\operatorname{Im}}

\newcommand{\Supp}{\operatorname{Supp}}

\newcommand{\Gal}{{\rm Gal}}

\newcommand{\Def}{\overset{{\rm def}}{=}}
\newcommand{\ab}{{\rm ab}}
\newcommand{\et}{\text{{\rm \'et}}}
\newcommand{\rmt}{\mathrm{t}}

\newcommand{\N}{\mathrm{N}}

\newcommand{\fppf}{{\rm fppf}}
\newcommand{\Diag}{\Delta}
\newcommand{\X}{\mathbb{X}}
\newcommand{\G}{\mathbb{G}}
\newcommand{\F}{\mathbb{F}}
\newcommand{\A}{\mathbb{A}}

\newcommand{\cosp}{\mathrm{cosp}}
\newcommand{\spe}{\mathrm{sp}}

\newcommand{\ols}{\overline{s}}
\newcommand{\olt}{\overline{t}}

\newcommand{\olx}{\overline{x}}

\newcommand{\olK}{\overline{K}}

\newcommand{\olR}{\overline{R}}
\newcommand{\M}{\mathcal{M}}

\newcommand{\oeta}{\overline{\eta}}

\newcommand{\scrO}{\mathscr{O}}

\newcommand{\scrC}{\mathscr{C}}

\newcommand{\scrD}{\mathscr{D}}

\newcommand{\calX}{\mathcal{X}}
\newcommand{\calY}{\mathcal{Y}}
\newcommand{\calZ}{\mathcal{Z}}
\newcommand{\calP}{\mathcal{P}}

\newcommand{\calF}{\mathcal{F}}
\newcommand{\lto}{\longrightarrow}
\newcommand{\cB}{\mathcal{B}}
\newcommand{\tame}{\mathrm{tame}}

\newcommand{\fX}{\mathfrak{X}}
\newcommand{\uAut}{\underline{\mathrm{Aut}}}
\newcommand{\uHom}{\underline{\mathrm{Hom}}}

\newcommand{\bfD}{\mathbf{D}}

\newcommand{\bfr}{\mathbf{r}}

\newcommand{\Pic}{\mathrm{Pic}}

\title[Tame fundamental group scheme]{The tame fundamental group schemes\\ of curves in positive characteristic}

\author[S. Otabe]{Shusuke Otabe}
\address{Graduate School of Engineering, Nagoya Institute of Technology, Nagoya, Aichi, Japan 
}
\email{shusuke.otabe@nitech.ac.jp}

\date{\today}


\thanks{}

\begin{document}

\begin{abstract}
The tame fundamental group scheme for an algebraic variety is the maximal linearly reductive quotient of Nori's fundamental group scheme. In this paper, we study the tame fundamental group schemes of smooth curves defined over algebraically closed fields of positive characteristic and develop the theory of cospecialization maps for them. As a result, we see that the tame fundamental group schemes heavily depend on the curves. We also see that numerical invariants of curves can be reconstructed from the tame fundamental group schemes.   
\end{abstract}

\maketitle

\thispagestyle{empty}

\section{Introduction}

Nori's \textit{fundamental group scheme} $\pi^{\N}(X)$~\cite[Chapter II]{No82} for a variety $X$ defined over a field $k$ is by definition a profinite $k$-group scheme classifying finite flat torsors of $X$. In the present paper, we will study its maximal \textit{linearly reductive} quotient $\pi^{\tame}(X)$~(cf.\ \S\ref{sec:tame fgs def}) in the case where $k=\overline{k}$ is an algebraically closed field of positive characteristic $p>0$. We call $\pi^{\tame}(X)$ the \textit{tame fundamental group scheme} following \cite[\S10]{BV15}, where the terminology `tame' comes from the notion of \textit{tame stacks} in the sense of \cite{AOV08}. As the group of $k$-valued points, we recover the maximal prime-to-$p$ quotient of Grothendieck's \textit{\'etale fundamental group}, i.e.\,$\pi^{\tame}(X)(k)\simeq\pi^{\et}_1(X)^{(p')}$. Therefore, the tame fundamental group scheme $\pi^{\tame}(X)$ is a group-scheme theoretic analogue of the prime-to-$p$ \'etale fundamental group $\pi^{\et}_1(X)^{(p')}$.

We will investigate the tame fundamental group schemes $\pi^{\tame}(U)$ of smooth curves $U$ defined over $k$. The structure of the \'etale part $\pi_1^{\et}(U)^{(p')}$ is well-understood by the experts (cf.\,\cite{Gr71}). If we denote by $g$ the genus of the smooth compactification $X=U^{\rm cpt}$ and by $n$ the cardinality $n=\#(X\setminus U)(k)$, then the isomorphism class of $\pi^{\et}_1(U)^{(p')}$ can be determined by the pair of integers $(g,n)$ and the profinite group $\pi_1^{\et}(U)^{(p')}$ is not dependent on the curve $U$. See \S\ref{sec:sp thm} for further details. In contrast, as main results of the present paper, we will see that the tame fundamental group scheme $\pi^{\tame}(U)$ heavily depends on the curve $U$.  As one of the main results, we will prove the following.

\begin{thm}(cf.\,Corollary \ref{cor:cosp2})\label{thm-int:main1}
Let $k_0=\overline{\F}_p$ be an algebraic closure of the prime field $\F_p$ of characteristic $p>0$ and set $S\Def\Spec k_0[[t]]=\{s,\eta\}$, where $s$ and $\eta$ are the closed point and the generic point of $S$, respectively. Let $X$ be a proper smooth relative $S$-curve of genus $g$ together with a relatively \'etale Cartier divisor $D$ on $X/S$ of degree $n$. We set $U\Def X\setminus\Supp(D)$. Suppose that $U$ is hyperbolic, i.e.\ $2-2g-n<0$. If $U_{\oeta}$ is not defined over $k_0$, then the tame fundamental group scheme $\pi^{\tame}(U_{\oeta})$ of $U_{\oeta}$ is not isomorphic to $\pi^{\tame}(U_s)\times_{k_0}k_0(\oeta)$. 
\end{thm}

We will also discuss reconstruction of numerical invariants (such as $g$, $n$) from the tame fundamental group schemes. For $i=1,2$, let $X_i$ be a proper smooth connected curve of genus $g_i$ and of $p$-rank $\gamma_i$ over $k$ and $S_i$ a finite set of closed points of $X_i$ with cardinality $n_i=\#S_i\ge 0$. We put $U_i\Def X_i\setminus S_i$. 
As the second main result, we will prove the following version of Tamagawa's theorem  \cite[Theorem 4.1]{Ta03}.

\begin{thm}(cf.\ Theorem \ref{thm:rcnst inv})\label{thm-int:main2}
Suppose that there exists an isomorphism of $k$-group schemes $\pi^{\tame}(U_1)\simeq\pi^{\tame}(U_2)$. Then we have $(g_1,n_1,\gamma_1)=(g_2,n_2,\gamma_2)$ unless $g_i=0$ and $n_i\le 1$ for $i=1,2$.
\end{thm}

We will explain the idea of the proof of Theorem \ref{thm-int:main1}. For the proof of the theorem, we will give an extension of the theory of \textit{specialization maps} for the \'etale fundamental groups to the tame fundamental group schemes. 
Let $S$ be a scheme. Let $X$ be a proper smooth relative $S$-curve of genus $g$ together with a relatively \'etale Cartier divisor $D$ on $X/S$ of degree $n$. Let $U\Def X\setminus \Supp(D)$ be the associated $S$-curve. We denote by $\pi^{\rmt}_1(U)$ the \textit{tame fundamental group} of $U$ in the sense of Grothendieck \cite[Expos\'e XIII, 2.1.3]{Gr71}. This is a profinite group classifying finite \'etale coverings of $U$ which are tamely ramified along $D$. By definition, the tame fundamental group $\pi_1^{\rmt}(U)$ is a quotient of the \'etale fundamental group $\pi^{\et}_1(U)$. In \cite[Expos\'e XIII, 2.10]{Gr71}, Grothendieck developed the theory of \textit{specialization maps} for the tame fundamental groups. The theory says that for any two geometric points $\ols,\olt$ of $S$, where $\ols$ is a specialization of $\olt$, there exists a surjective continuous homomorphism, called the \textit{specialization map},
\begin{equation*}
\spe^{\rmt}\colon\pi^{\rmt}_1(U_{\olt})\twoheadrightarrow\pi^{\rmt}_1(U_{\ols})
\end{equation*}
which is canonical up to conjugation (see also \S\ref{sec:sp thm}). The theory of the specialization maps has many applications to analysis of the tame fundamental groups of curves defined over fields of positive characteristic $p>0$. Indeed, Grothendieck himself used it to prove the finite generatedness of the tame fundamental groups and to give an explicit description of the maximal pro-prime-to-$p$ quotient $\pi^{\et}_1(U)^{(p')}$ of the \'etale fundamental group (cf.\ Corollary \ref{cor:pi^t fin gen}). Other uses appear in anabelian geometry (cf.\ \cite{Ra02}\cite{PS03}\cite{Ta04}\cite{Sa05}; see also \S\ref{sec:Tamagawa's thm}), where the specialization map plays a crucial role to prove non-constancy results for the tame fundamental group on the moduli space $\M_{g,[n],\F_p}$ of $n$-pointed genus $g$ curves.

To obtain an analogous theory for the tame fundamental group schemes, we will  consider a smaller quotient $\pi^{\tame}(U)^{\scrD}$, which is the maximal pro-$\scrD$ quotient of the tame fundamental group scheme, where $\scrD$ is the class of finite linearly reductive group schemes $G$ whose connected part $G^0$ is \textit{elementary}, i.e.\ $G^0\simeq\mu_p^{s}$ for some integer $s\ge 0$ (see Definition \ref{def:category D}). As a result, we will establish the following theorem.

\begin{thm}(cf.\,Proposition \ref{prop:cosp}, Theorem \ref{thm:cosp})\label{thm-int:key}
Let $k$ be an algebraically closed field of characteristic $p>0$ and set $S\Def\Spec k[[t]]=\{s,\eta\}$, where $s$ and $\eta$ are the closed point and the generic point of $S$ respectively. Let $X$ be a proper smooth relative $S$-curve of genus $g$ together with a relatively \'etale Cartier divisor $D$ on $X/S$ of degree $n$. We set $U\Def X\setminus\Supp(D)$. 
\begin{enumerate}
\renewcommand{\labelenumi}{(\arabic{enumi})}
\item There exists a canonical homomorphism of affine $k(\oeta)$-group schemes
\begin{equation*}
\cosp^{\scrD}\colon \pi^{\tame}(U_s)^{\scrD}\times_{k}k(\oeta)\to\pi^{\tame}(U_{\oeta})^{\scrD}
\end{equation*}
up to conjugation, which we call the \textit{cospecialization map}, such that the following conditions are satisfied.
\begin{enumerate}
\renewcommand{\labelenumii}{(\roman{enumii})}
\item The map $\cosp^{\scrD}$ is injective.

\item By taking the groups of $k(\oeta)$-valued points, the map $\cosp^{\scrD}$ induces an  isomorphism between the maximal pro-prime-to-$p$ quotients of the \'etale fundamental groups~(cf.\ (\ref{eq:sp p'})) $\pi^{\et}_1(U_s)^{(p')}\simeq\pi_1^{\et}(U_{\oeta})^{(p')}$.
\end{enumerate} 
\item The following are equivalent.
\begin{enumerate}
\renewcommand{\labelenumi}{(\alph{enumi})}
\item The cospecialization map $\cosp^{\scrD}$ is isomorphism. 
\item For any $G\in\scrD$, we have $\#\Hom(\pi^{\tame}(U_s,x_s)^{\scrD},G)=\#\Hom(\pi^{\tame}(U_{\oeta},x_{\oeta})^{\scrD},G_{k(\oeta)})$.
\item There exists an isomorphism of $k$-group schemes $\pi^{\tame}(U_s,x_s)^{\scrD}\times_kk(\oeta)\simeq\pi^{\tame}(U_{\oeta},x_{\oeta})^{\scrD}$.
\end{enumerate}
\end{enumerate}  
\end{thm}

\begin{rem}
\begin{enumerate}
\renewcommand{\labelenumi}{(\arabic{enumi})}
\item Theorem \ref{thm-int:key}(1) is a variant of Grothendieck's specialization theorem (cf.\,\cite[Expos\'e XIII, 2.10]{Gr71}; see also \S\ref{sec:sp thm}). The surjectivity of the specialization map $\spe^{\rmt}\colon\pi^{\rmt}_1(U_{\olt})\twoheadrightarrow\pi^{\rmt}_1(U_{\ols})$ is replaced by the injectivity of the cospecialization map. Theorem \ref{thm-int:key}(2) is a key ingredient for the proof of Theorem \ref{thm-int:main1}. 
Other ingredients of the proof of Theorem \ref{thm-int:main1} are structural results discussed in \S\ref{sec:str tame fg} and Tamagawa's specialization theorem for proper smooth curves (cf.\,\cite[Theorem 6.1]{Ta04}; see also Theorem \ref{thm:Tamagawa C}). Theorem \ref{thm-int:main1} can be considered as a version of Tamagawa's specialization theorem (cf.\,\cite[Theorem 8.1]{Ta04}; see Theorem \ref{thm:Tamagawa}). 

\item One can also see that if $U_{\oeta}$ is isomorphic to the trivial deformation $U_{s}\times_{k}k(\oeta)$ of the special fiber $U_{s}$, then the cospecialization map becomes an isomorphism $\cosp^{\scrD}\colon \pi^{\tame}(U_{s})^{\scrD}\times_{k}k(\oeta)\xrightarrow{\simeq}\pi^{\tame}(U_{\oeta})^{\scrD}$, which can be deduced from the base change theorem for the tame fundamental group schemes of curves (cf.\,Proposition \ref{prop:bc tame fgs}). 
On the other hand, Theorem \ref{thm-int:main1} implies that the tame fundamental group scheme is not constant on the moduli space $\M_{g,[n],k_0}$ of $n$-pointed genus $g$ curves over $k_0$ under the assumption that $2-2g-n<0$. However, this result is much weaker than the non-constancy result for the tame fundamental group due to Tamagawa~(cf.\ \cite[Theorem 8.6]{Ta04}), which asserts that the tame fundamental group is not constant on the set $\M_{g,[n]}(k_0)$ of $k_0$-valued points.   

\item As an application of \cite[Theorem 4.1]{Ta03}, Tamagawa provided a group-theoretic characterization of inertia groups of the tame fundamental group $\pi_1^{\rmt}(U)$ for affine hyperbolic curves $U/k$ (cf.\,\cite[Theorem 5.2]{Ta03}). It seems natural to seek an analogous use of Theorem \ref{thm-int:main2} for the tame fundamental group scheme $\pi^{\tame}(U)$. However, the author has no idea to reconstruct \textit{inertia subgroup schemes} in a purely group-scheme theoretic manner at present.

\item The starting point of this work was the expectation that results known for prime-to-$p$ Galois coverings could be extended to finite linearly reductive torsors in some form. However, lifting problems for $\mu_p$-torsors already indicate that the situation is quite different. When $U=X$, the Serre duality implies that there exists a natural isomorphism $H^{1}_{\fppf}(X_{\olt},\mu_p)\simeq H_{\et}^1(X_{\olt},\Z/p\Z)^{\vee}$ for $t\in\{s,\eta\}$,  which suggests that the classification of $\mu_p$-torsors is closely related to the classification of $\Z/p\Z$-torsors. However, the geometries of them are different. Indeed, it is known that any $\mu_p$-torsor over $X_{\oeta}$ can be extended to a $\mu_p$-torsor over $X\times_S S'$ for some finite extension $S'\to S$. On the other hand, the same does not hold for $\Z/p\Z$-coverings. The non-existence of $\Z/p\Z$-models over $X$ gives an obstruction for the injectivity of the specialization map $\spe\colon \pi^{\et}_1(X_{\oeta})\twoheadrightarrow\pi_1^{\et}(X_{s})$. Instead, any $\Z/p\Z$-coverings over $X_{s}$ can be uniquely lifted to a $\Z/p\Z$-covering over $X$, but deformations of a $\mu_p$-torsor of $X_{\ols}$ are far from unique. The non-uniqueness of deformations of $\mu_p$-torsors gives an obstruction for the surjectivity of the cospecialization map $\cosp^{\scrD}\colon\pi^{\tame}(X_{s})^{\scrD}\times_{k}{k(\oeta)}\hookrightarrow\pi^{\tame}(X_{\oeta})^{\scrD}$.
\end{enumerate}
\end{rem}

We end this introduction section with the organization of the present paper. 
In \S\ref{sec:pre}, we recall some basic notions which we freely use in the present paper. In \S\ref{sec:Nori fg}, we recall the definition of the \textit{Nori fundamental gerbe} following \cite{BV15} (see also \cite[Chapter II]{No82}). In \S\ref{sec:root stack}, we recall the definition of \textit{root stacks} in the sense of \cite{AGV08}. In \S\ref{sec:etale fg}, we recall Grothendieck's construction of the specialization map for the tame fundamental group (cf.\,\S\ref{sec:sp thm}). We also recall the theorem of Tamagawa (cf.\,\S\ref{sec:Tamagawa's thm}). In the final subsection (cf.\ \S\ref{sec:fq et fg}), we recall a description of certain finite quotients of the \'etale fundamental groups of proper smooth curves. In \S\ref{sec:tame fgs}, we recall the definition of the {tame fundamental group scheme} (cf.\,\S\ref{sec:lin red}, \S\ref{sec:tame fgs def}) and discuss several structural results (cf.\,\S\ref{sec:str tame fg}, \S\ref{sec:bc tame fgs}).

In \S\ref{sec:main}, we prove the main results. First we recall a description of finite quotients of the tame fundamental group scheme in \S\ref{sec:fq tame fgs}. Next we settle a certain lifting problem for finite linearly reductive torsors of curves in \S\ref{sec:LP}, which is crucial for the condition (i) in Theorem \ref{thm-int:key}(1). In the  subsection \S\ref{sec:cosp}, we construct the cospecialization map and prove Theorems \ref{thm-int:main1} and \ref{thm-int:key}. In the final subsection \S\ref{sec:rmk inv}, we prove Theorem \ref{thm-int:main2}.

\begin{ack}
The author would like to express his gratitude to Prof.\,Takao Yamazaki, Prof.\,Lei Zhang and Prof.\,Yu Yang for fruitful discussions and for their advice. Especially, a Zoom meeting with Prof.\,Yang in January 2021 led the author to prove Theorem \ref{thm-int:main2}. A part of this work was done during the author's visit to Freie Universit\"at Berlin in 2018. He would like to thank Prof.\,H\'el\`ene Esnault and all  the members of her workgroup for their hospitality. The author would like to thank the referees for reading this manuscript and for giving helpful comments. This work was supported by JSPS KAKENHI Grant Numbers JP16J02171, JP19J00366,  JP21K20334, JP24K16894. The stay in Berlin was supported by JSPS Overseas Challenge Program for Young Researchers. 
\end{ack}


\section{Preliminaries}\label{sec:pre}

For a field $k$, we denote by $\Vecf_k$ the category of finite dimensional vector spaces over $k$. For an affine $k$-group scheme $G$, we denote by $\Rep(G)$ the category of finite dimensional left $k$-linear representations of $G$. 
For an algebraic stack $\calX$ over a scheme $S$, we denote by $\Qcoh(\calX)$ (respectively $\Vect(\calX)$) the category of quasi-coherent sheaves on $\calX$ (respectively the category of vector bundles on $\calX$). For a field $k$, we have $\Vect(\Spec k)=\Vecf_k$. If $\calX$ is Noetherian, we also consider the category $\Coh(\calX)$ of coherent sheaves on $\calX$. If $\calX=\cB_S G$ is the classifying stack of an affine flat and finitely presented $S$-group scheme $G$, then  $\Qcoh(\cB_S G)$ is nothing but the category of \textit{$G$-equivariant} quasi-coherent sheaves on $S$. Namely, it is the category of quasi-coherent sheaves $\calF$ on $S$ endowed with an action of $G$~(cf.~\cite[\S2.1]{AOV08}). In the case where $S=\Spec k$ is the spectrum of a field $k$, all the three categories $\Coh(\cB G)$, $\Vect(\cB G)$ and $\Rep(G)$ are canonically equivalent to each other.

\subsection{The Nori fundamental gerbe}\label{sec:Nori fg}

\begin{para}\label{para:fin stack}
Let $k$ be a field. A \textit{finite stack} over $k$ is an algebraic stack $\Gamma$ over $k$ which has finite flat diagonal and admits a flat surjective morphism $U\to\Gamma$ for some finite $k$-scheme $U$~(cf.~\cite[Definition 4.1]{BV15}). A \textit{finite gerbe} over $k$ is a finite stack over $k$ which is a gerbe in the fppf topology. A finite stack $\Gamma$ is a finite gerbe if and only if it is geometrically connected and geometrically reduced~(cf.~\cite[Proposition 4.3]{BV15}). A \textit{profinite gerbe} over $k$ is 
a projective limit of finite gerbes over $k$~(cf.~\cite[Definition 4.6]{BV15}).  
\end{para}

\begin{para}\label{para:Nori FG}
Let $\calX$ be an algebraic stack of finite type over $k$. Suppose that $\calX$ is \textit{inflexible} in the sense of \cite[Definition 5.3]{BV15}. For example, if $\calX$ is geometrically connected and geometrically reduced, then it is inflexible~(cf.\ \cite[Proposition 5.5(b)]{BV15}). Then, there exists a profinite gerbe $\Pi$ over $k$ together with a morphism $\calX\to\Pi$ such that, for any finite stack $\Gamma$ over $k$, the induced functor
\begin{equation*}
\Hom_k(\Pi,\Gamma)\to\Hom_k(\calX,\Gamma)
\end{equation*} 
is an equivalence of categories~(cf.~\cite[Theorem 5.7]{BV15}). Such a gerbe $\Pi$ is unique up to unique isomorphism, so we denote it by $\Pi_{\calX/k}^{\N}$, and call it the \textit{Nori fundamental gerbe} for $\calX$ ~(cf.~\cite{BV15}). If $G$ is a finite $k$-group scheme, then the associated classifying stack $\cB G$ is a finite gerbe and there exists a natural bijection
\begin{equation*}
\Hom_k(\Pi_{\calX/k}^{\N},\cB_kG)\xrightarrow{\simeq}\Hom_k(\calX,\cB_kG)=H_{\fppf}^1(\calX,G).
\end{equation*}
\end{para}

\begin{para}\label{para:Nori fgs}
If $\calX$ admits a $k$-rational point $x\in \calX(k)$, then the composition $\Spec k\xrightarrow{x}\calX\to\Pi_{\calX/k}^{\N}$ defines a section $\xi\in\Pi_{\calX/k}^{\N}(k)$. We denote by $\pi^{\N}(\calX,x)$ the automorphism group scheme $\underline{\mathrm{Aut}}_k(\xi)$, i.e.\ $\pi^{\N}(\calX,x)\Def\uAut_k(\xi)$. Let $X_x^{\N}\to X$ be the fpqc $\pi^{\N}(\calX,x)$-torsor associated with the morphism $\calX\to\Pi^{\N}_{\calX/k}\xleftarrow{\simeq}\cB_k\pi^{\N}(\calX,x)$. By definition, it admits a unique $k$-rational point $x^{\N}\in X_x^{\N}(k)$ above $x$. The resulting triple $(X^{\N}_x,\pi^{\N}(X,x),x^{\N})$ then recovers Nori's construction of the \textit{fundamental group scheme} of $(X,x)$ in \cite[Chapter II]{No82}. Namely, for any finite $k$-group scheme $G$, the set of homomorphisms $\Hom(\pi^{\N}(\calX,x),G)$ is naturally bijective onto the set of isomorphism classes of pointed $G$-torsors $(P,p)\to(\calX,x)$. More precisely, for each homomorphism $\phi\colon\pi^{\N}(\calX,x)\to G$ into a finite $k$-group scheme, the corresponding pointed $G$-torsor is given by $(P,p)\Def (X_x^{\N},x^{\N})\times^{\pi^{\N}(\calX,x)_{\phi}}G$.
\end{para}

\begin{para}\label{para:Nori-reduced}
Let $\calX$ be an inflexible algebraic stack of finite type over $k$. A morphism $\calX\to\Gamma$ into a finite gerbe $\Gamma$ is said to be \textit{Nori-reduced}~(cf.~\cite[Definition 5.10]{BV15}) if for any factorization $\calX\to
\Gamma'\to\Gamma$ where $\Gamma'$ is a finite gerbe and $\Gamma'\to\Gamma$ is faithful, then $\Gamma'\to\Gamma$ is an isomorphism. According to \cite[Lemma 5.12]{BV15}, for any morphism $\calX\to\Gamma$ into a finite gerbe, there exists a unique factorization $\calX\to\Delta\to\Gamma$, where $\Delta$ is a finite gerbe,  $\calX\to\Delta$ is Nori-reduced and $\Delta\to\Gamma$ is representable. A $G$-torsor $\calP\to \calX$ is said to be \textit{Nori-reduced} if the morphism $\calX\to\cB G$ is Nori-reduced.
\end{para}

\begin{para}\label{para:EFin}
Under the assumption that $\calX$ is proper over $k$, the Nori fundamental gerbe $\Pi^{\N}_{\calX/k}$ has a tannakian interpretation in terms of vector bundles on $\calX$. Indeed, the pullback functor of the morphism $\calX\to\Pi^{\N}_{\calX/k}$ induces a fully faithful tensor functor
\begin{equation*}\label{eq:EFin}
\Vect(\Pi^{\N}_{\calX/k})\to\Vect(\calX)
\end{equation*}
whose essential image is the tannakian category 
 $\EFin(\calX)$ of \textit{essentially finite bundles} on $\calX$~(cf.\ \cite[\S7]{BV15}), i.e.\
\begin{equation}\label{eq:EFin}
\Vect(\Pi^{\N}_{\calX/k})\xrightarrow{\simeq}\EFin(\calX)\subset\Vect(\calX).
\end{equation}
\end{para}

\subsection{Root stacks}\label{sec:root stack}

In this subsection, we will recall the definition of \textit{root stacks}~(cf.\ \cite[Appendix B.2]{AGV08}). 
Let $S$ be a scheme and $X$ a Noetherian $S$-scheme. Let $\bfD=(D_i)_{i=1}^n$ be an $n$-tuple of reduced irreducible relative effective Cartier divisors $D_i$ on $X/S$. For each $i$, let $\scrO_X(D_i)$ be the line bundle associated with $D_i$ and $s_{D_i}\in\Gamma(X,\scrO_X(D_i))$ a canonical section. Then the pair $(\scrO_X(D_i),s_{D_i})$ gives rise to a morphism $\phi_{i}\colon X\to[\A^1_S/\G_{m,S}]$ into the quotient stack $[\A^1_{S}/\G_{m,S}]$ of the affine line $\A^1_S$ with respect to the natural action by the multiplicative group scheme $\G_{m,S}$. By taking the fiber product over $S$, we get the $S$-morphism 
\begin{equation*}
\phi_{(X,\bfD)}\Def(\phi_i)_{i=1}^n\colon X\to\prod_{i=1}^n[\A^1_S/\G_{m,S}]=[\A^n_{S}/\G_{m,S}^n].
\end{equation*}
Here each morphism $\phi_{i}$ does not depend on the choice of the canonical section $s_{D_i}$ and the morphism $\phi_{(X,\bfD)}$ is natural with respect to the pair $(X,\bfD)$. 

\begin{para}\label{para:def of root stack}
Now for any $n$-tuple of positive integers $\bfr=(r_i)_{i=1}^n$, by taking the $2$-fiber product of the $\bfr$-th power map $\theta_{\bfr}\colon[\A_S^n/\G_{m,S}^n]\to[\A^n_{S}/\G_{m,S}^n]$, we obtain the \textit{root stack} $\sqrt[\bfr]{\bfD/X}$ associated with the data $(X,\bfD,\bfr)$,
\begin{equation*}
\begin{xy}
\xymatrix{\ar@{}[rd]|{\square}
\sqrt[\bfr]{\bfD/X}\ar[r]\ar[d]_{\pi}&[\A^n/\G_{m}^n]\ar[d]^{\theta_{\bfr}}\\
X\ar[r]&[\A^n/\G_{m}^n].
}
\end{xy}
\end{equation*}

The root stack $\sqrt[\bfr]{\bfD/X}$ is a \textit{tame stack} over $S$ in the sense of \cite{AOV08}, where the natural projection map $\pi\colon\sqrt[\bfr]{\bfD/X}\to X$ gives the coarse moduli space map. In particular, the push-forward functor of the categories of quasi-coherent sheaves $\pi_*\colon\Qcoh(\sqrt[\bfr]{\bfD/X})\to\Qcoh(X)$ is exact.

As the coarse moduli space map $\pi$ is proper by definition, if $X$ is proper over $S$, then so is the root stack $\sqrt[\bfr]{\bfD/X}$. Note that the map $\pi\colon\sqrt[\bfr]{\bfD/X}\to X$ is an isomorphism over the open subscheme $U\Def X\setminus D$, where $D\Def\bigcup_{i=1}^n {\rm Supp}(D_i)$ and hence we have a natural open embedding $U\hookrightarrow\sqrt[\bfr]{\bfD/X}$. 
\end{para}

\begin{para}\label{para:res gerbe}
On the other hand, the local picture around a stacky point can be described as follows. Namely, for any closed point $x\in|\sqrt[\bfr]{\bfD/X}|_0=|X|_0$ with $x\not\in |U|_0$, there exists a closed immersion $\cB\mu_{\bfr_x,k(x)}\hookrightarrow\sqrt[\bfr]{\bfD/X}$, where $\bfr_x=(r_i)_{x\in D_i}$ which fits into the following commutative diagram
\begin{equation*}
\begin{xy}
\xymatrix{
&\cB\mu_{\bfr_x,k(x)}\ar@{^{(}->}[r]\ar[d]&\sqrt[\bfr]{\bfD/X}\ar[d]^{\pi}\\
&\Spec k(x)\ar[r]^{~~~x}&X.
}
\end{xy}
\end{equation*}
The closed immersion $\cB\mu_{\bfr_x,k(x)}\hookrightarrow\sqrt[\bfr]{\bfD/X}$ is nothing but the \textit{residual gerbe} at the point $x$ in the sense of \cite[Definition 06MU]{stack}. In particular, the residual gerbes of the root stacks are always neutral gerbes.  
\end{para}


\section{The specialization map for the tame fundamental group}\label{sec:etale fg}

\subsection{Grothendieck's specialization theorem}\label{sec:sp thm}

In this subsection, we will recall Grothendieck's specialization theorem for the tame fundamental group~(cf.\ \cite[Expos\'e XIII, 2.10]{Gr71}). Let $S$ be a scheme. Let $X$ be a proper smooth relative $S$-curve of genus $g$ with geometrically connected fibers and $D$ a relatively \'etale Cartier divisor on $X/S$ of degree $n$. We set $U\Def X\setminus \Supp(D)$. For each geometric point $\olx\to U$, we denote by $\pi_1^{\rmt}(U,\olx)$ the \textit{tame fundamental group} of $U$ with respect to $\olx$ (cf.\ \cite[Expos\'e XIII, 2.1.3]{Gr71}). Note that if $D=\emptyset$ or $S$ is the spectrum of a field of characteristic $0$, then the tame fundamental group $\pi^{\rmt}_1(U,\olx)$ coincides with the \'etale fundamental group $\pi^{\et}_1(U,\olx)$. Let $\ols,\olt$ be geometric points of $S$ such that $\ols$ is a specialization of $\olt$. Let $\widetilde{S}$ be the strict henselization of $S$ at $\ols$. Then we get the commutative diagram
\begin{equation*}
\begin{xy}
\xymatrix{
U_{\olt}\ar[r]\ar[d]&\widetilde{U}\ar[d]&U_{\ols}\ar[d]\ar[l]\\
\olt\ar[r]&\widetilde{S}&\ols,\ar[l]
}
\end{xy}
\end{equation*}
with the Cartesian squares, where $\widetilde{U}\Def U\times_S\widetilde{S}$.

If we choose the geometric points ${\olx}_1$ and $\olx_2$ of $U_{\olt}$ and $U_{\ols}$ respectively, then we have canonical homomorphisms between the tame fundamental groups
\begin{equation*}
\phi_1\colon\pi^{\rmt}_1(U_{\olt},\olx_1)\twoheadrightarrow\pi_1^{\rmt}(\widetilde{U},\olx_1)\quad\text{and}\quad\phi_2\colon\pi^{\rmt}_1(U_{\ols},\olx_2)\xrightarrow{\simeq}\pi_1^{\rmt}(\widetilde{U},\olx_2),
\end{equation*}
where the first homomorphism $\phi_1$ is surjective and the second one $\phi_2$ is an isomorphism~(cf.\ \cite[Expos\'e XIII, 2.10]{Gr71}). Therefore, if we take any path $\phi_{12}\colon\pi^{\rmt}_1(\widetilde{U},\olx_1)\xrightarrow{\simeq}\pi_1^{\rmt}(\widetilde{U},\olx_2)$ from $\olx_1$ to $\olx_2$, then we get a surjective continuous homomorphism of profinite groups
\begin{equation}\label{eq:sp^t}
\spe^{\rmt}\Def\phi_2^{-1}\phi_{12}\phi_1\colon\pi^{\rmt}_1(U_{\olt},\olx_1)\twoheadrightarrow\pi^{\rmt}_1(U_{\ols},\olx_2),
\end{equation}
which we call the \textit{specialization map} for the tame fundamental group. By definition, the map $\spe^{\rmt}$ is canonically determined by the curve $U$ together with the geometric points $\ols,\olt$ of $S$ up to inner automorphism of $\pi^{\rmt}_1(U_{\ols},\olx_2)$. Thus, we will often write
\begin{equation*}
\spe^{\rmt}\colon\pi_1^{\rmt}(U_{\olt})\to\pi^{\rmt}_1(U_{\ols})
\end{equation*}
without mentioning the base points. 

In the case where $D=\emptyset$, i.e.\ $X=U$, we get the specialization map for the \'etale fundamental group 
\begin{equation}\label{eq:sp}
\spe=\spe^{\rmt}\colon\pi^{\et}_1(X_{\olt},\olx_1)\twoheadrightarrow\pi^{\et}_1(X_{\ols},\olx_2).
\end{equation}

Moreover, if the point $\ols$ has residue characteristic $p>0$, then the specialization map (\ref{eq:sp^t})   induces an isomorphism between the maximal pro-prime-to-$p$ quotients of the \'etale fundamental groups
\begin{equation}\label{eq:sp p'}
\spe^{(p')}\colon\pi^{\et}_1(U_{\olt},\olx_1)^{(p')}\xrightarrow{\simeq}\pi^{\et}_1(U_{\ols},\olx_2)^{(p')}.
\end{equation}

These are the contents of 
Grothendieck's specialization theorem. 
As an application, we have the following consequence.

\begin{cor}(cf.\ \cite[Expos\'e XIII, Corollaire 2.12]{Gr71})\label{cor:pi^t fin gen}
Let $k$ be an algebraically closed field of characteristic $p>0$. Let $X_0$ be a proper smooth connected curve of genus $g$ over $k$ together with a finite (possibly empty) set $D_0$ of closed points of $X_0$. We set $n\Def\# S$. Let $X$ be a proper smooth relative curve of genus $g$ over $S\Def\Spec W(k)$ and $D$ a set of $S$-valued points of $X$ such that $(X,D)\times_{S}\Spec k\simeq (X_0,D_0)$. We set $U\Def X\setminus D$ and $U_0\Def X_0\setminus D_0$. Let $\eta\in S$ be the generic point of $S$. Then there exists a surjective continuous homomorphism
\begin{equation*}
\pi_1^{\et}(U_{\oeta})\twoheadrightarrow\pi_1^{\rmt}(U_0)
\end{equation*}
which induces an isomorphism between the maximal pro-prime-to-$p$ quotients
\begin{equation*}
\pi_1^{\et}(U_{\oeta})^{(p')}\xrightarrow{\simeq}\pi_1^{\et}(U_0)^{(p')}.
\end{equation*}
In particular, the tame fundamental group $\pi_1^{\rmt}(U_0)$ of the curve $U_0$ is topologically finitely generated and we have an isomorphism of pro-prime-to-$p$ groups
\begin{equation*}
\pi^{\et}_1(U_0)^{(p')}\simeq\Pi_{g,n}^{(p')},
\end{equation*} 
where we define
\begin{equation}\label{eq:Pi g,n}
\Pi_{g,n}\Def
\biggl\langle 
\begin{aligned}
a_1,b_1,&\dots,a_g,b_g,\\
\delta_1,&\dots,\delta_n
\end{aligned}
~\biggl |~\prod_{i=1}^g[a_i,b_i]\delta_1\cdots\delta_n=1\biggl\rangle.
\end{equation}
\end{cor}

\subsection{Tamagawa's specialization theorem}\label{sec:Tamagawa's thm}

In this subsection, we will recall use of the specialization maps in anabelian geometry for hyperbolic curves over algebraically closed fields of positive characteristic~(cf.\ \cite{Ra02}\cite{PS03}\cite{Ta04}\cite{Sa05}). Especially, we will recall Tamagawa's specialization theorems (see Theorems \ref{thm:Tamagawa} and \ref{thm:Tamagawa C} below).

Let $k_0\Def\overline{\F}_p$. Let $S$ be an $\F_p$-scheme and $U$ a smooth relative $S$-curve as in \S\ref{sec:sp thm}. Recall that the curve $U$ is said to be \textit{hyperbolic} if $2-2g-n<0$. Moreover, a curve over a field containing $k_0$ is said to be \textit{constant} if it is defined over $k_0$.

\begin{thm}(cf.\ \cite[Theorem 8.1]{Ta04})\label{thm:Tamagawa}
Let $U$ be a hyperbolic $S$-curve and $\ols,\olt\in S$ two geometric points of $S$ such that $\ols$ is a specialization of $\olt$. Suppose that $U_{\ols}$ is constant and that $U_{\olt}$ is not constant. Then the specialization map $\spe^{\rmt}\colon\pi_1^{\rmt}(U_{\olt})\to\pi_1^{\rmt}(U_{\ols})$ for the tame fundamental group~(cf.\ (\ref{eq:sp^t})) is not an isomorphism. 
\end{thm}

As a corollary of the theorem, Tamagawa also obtained a non-constancy result of the tame fundamental group on the moduli space $\M_{g,[n],\F_p}$ of $n$-pointed genus $g$ curves~(cf.\ \cite[Theorems 8.3 and 8.6]{Ta04}).

For later use, let us reformulate the theorem for proper curves as follows.

\begin{definition}\label{def:category C}
Let $\scrC$ be the category of finite groups which has an elementary abelian normal $p$-Sylow subgroup. For any profinite group $\Pi$, we define the quotient $\Pi^{\scrC}$ to be
\begin{equation}\label{eq:Pi^C}
\Pi^{\scrC}\Def\varprojlim_{N\in\scrC(\Pi)}\Pi/N,
\end{equation}
where $\scrC(\Pi)$ is the set of open normal subgroups $N$ of $\Pi$ such that $\Pi/N\in\scrC$. 
\end{definition}

Note that for any $G\in\scrC$ with $P\triangleleft G$ the unique $p$-Sylow subgroup, the short exact sequence
\begin{equation*}
1\to P\to G\to G/P\to 1
\end{equation*}
splits, i.e.\ $G\simeq (G/P)\ltimes P$.

\begin{thm}(cf.\ \cite[Theorem 6.1 and Remark 6.3]{Ta04}\cite[Proposition 2.2.4(2)]{Ra02})\label{thm:Tamagawa C}
Let $S=\Spec k_0[[t]]=\{s,\eta\}$, where $s$ (respectively $\eta$) is the closed point (respectively the generic point) of $S$. Let $X$ be a proper smooth relative $S$-curve of genus $g\ge 2$ with geometrically connected fibers. Suppose that $X_{\oeta}$ is not constant. Then the map $\spe^{\scrC}\colon\pi_1^{\et}(X_{\olt})^{\scrC}\twoheadrightarrow\pi^{\et}_1(X_{\ols})^{\scrC}$ induced by the specialization map for the \'etale fundamental group (cf.\ (\ref{eq:sp})) is not an isomorphism.
\end{thm}

\begin{rem}\label{rem:Tamagawa C}
In fact, the quotient considered in \cite[Proposition 2.2.4(2)]{Ra02} or \cite[Remark 6.3]{Ta04} is not our $\pi^{\et}_1(X)^{\scrC}$ but a larger one $\pi_1^{(p,p')}(X)$~(cf.\ \cite[\S2, p.\ 343]{Ra02}). Precisely, if we set $P\Def \Ker(\pi_1^{(p,p')}(X)\twoheadrightarrow\pi_1^{\et}(X)^{(p')})$, we have
\[
\pi^{\et}_1(X)^{\scrC}=\pi_1^{(p,p')}(X)/[P,P]P^p.
\]
The latter group is the maximal pro-$\scrC'$ quotient of $\pi_1^{\et}(X)$, i.e.\ $\pi^{(p,p')}_1(X)\Def\pi_1^{\et}(X)^{\scrC'}$, where $\scrC'$ is the class of finite groups which has a normal $p$-Sylow subgroup. Note that $\scrC\subset\scrC'$ and that each group $G\in\scrC'$ is isomorphic to a semi-direct product $G\simeq H\ltimes P$, where $H$ has prime-to-$p$ order and $P$ is a $p$-group~(cf.\ \cite[1.1]{PS00}).

However, the isomorphism class of the profinite group $\pi^{(p,p')}_1(X)=\pi^{\et}_1(X)^{\scrC'}$ can be completely determined by the smaller quotient $\pi_1^{\et}(X)^{\scrC}$. Indeed, as $\pi_1^{\et}(X)^{\scrC'}$ is topologically finitely generated~(cf.\ Corollary \ref{cor:pi^t fin gen}), the isomorphism class of the profinite group $\pi_1^{\et}(X)^{\scrC'}$ can be determined by the set $\pi_A^{\et}(X)^{\scrC'}$ of finite quotients of $\pi_1^{\et}(X)^{\scrC'}$~(cf.\ \cite[Lemma 8.4]{Ta04}; see also \S\ref{sec:fq et fg}). However, by the result of Pacheco--Stevenson~\cite[Theorem 1.3]{PS00} (see also \cite[Proposition 2.5]{Bo01}), the set of finite quotients $\pi^{\et}_A(X)^{\scrC'}$ can be completely determined by the one $\pi_A^{\et}(X)^{\scrC}$ of the maximal pro-$\scrC$ quotient $\pi^{\et}_1(X)^{\scrC}$. Therefore, Theorem \ref{thm:Tamagawa C} is a valid reinterpretation of Tamagawa's theorem~\cite[Theorem 6.1]{Ta04}.
\end{rem}

\subsection{Finite quotients of the \'etale fundamental groups of proper curves}\label{sec:fq et fg}

In the whole of this subsection, we fix an algebraically closed field $k$ of characteristic $p>0$.

Let $X$ be a proper smooth connected curve of genus $g$ over $k$. As the \'etale fundamental group $\pi^{\et}_1(X)$ is topologically finitely generated~(cf.\ Corollary \ref{cor:pi^t fin gen}), the isomorphism class of the profinite group $\pi^{\et}_1(X)$ can be completely determined by the set $\pi^{\et}_A(X)$ of isomorphism classes of finite quotients of $\pi^{\et}_1(X)$~(cf.\ \cite[Lemma 8.4]{Ta04}). 

We denote by $\pi^{\et}_A(X)^{\scrC}$ is the set of isomorphism classes of finite quotients of the profinite group $\pi_1^{\et}(X)^{\scrC}$~(cf.\ (\ref{eq:Pi^C})). As Theorem \ref{thm:Tamagawa C} suggests, the set $\pi^{\et}_A(X)^{\scrC}$ of finite quotients is already complicated. This is in fact caused by the complexity of the $p$-\textit{rank} of prime-to-$p$ Galois coverings over $X$~(cf.\ \cite{Na83}\cite{Ruck}\cite{Pa95}\cite{Bo01}).

\begin{definition}(cf.\ \cite[\S1]{Pa95}\cite{Bo01})\label{def:p-rank}
The $p$-\textit{rank} $\gamma_A$ of an abelian variety $A$ over $k$ is defined to be $\gamma_A\Def\Dim_{\F_p}A[p](k)$. 
The $p$-\textit{rank} $\gamma_X$ of a proper smooth connected curve $X$ over $k$ is defined as the $p$-rank $\gamma_J$ of the Jacobian variety $J=J_X$ of $X$. The $p$-rank $\gamma_X$ is also called the \textit{Hasse--Witt invariant} for $X$. 
\end{definition}

Let $X$ be a proper smooth connected curve of genus $g$ over $k$. The Kummer theory in the fppf topology gives an isomorphism
\begin{equation}\label{eq:p-rank mu_p}
H^1_{\fppf}(X,\mu_p)\xrightarrow{\simeq}J[p](k).
\end{equation}
Moreover, according to \cite[III, Proposition 4.14]{Mi80}, we also have a natural isomorphism $H^1_{\fppf}(X,\mu_p)\simeq H^0(X,\Omega^1_X)^{\mathcal{C}}$, where $\mathcal{C}$ is the Cartier operator and the right hand side is the $\F_p$-subspace of $H^0(X,\Omega^1_X)$ consisting of regular differential forms $\omega\in H^0(X,\Omega^1_X)$ satisfying $\mathcal{C}(\omega)=\omega$. On the other hand, by the Artin--Schreier theory, we have an isomorphism
of $\F_p$-vector spaces $H^1_{\et}(X,\Z/p\Z)\simeq H^1(X,\scrO_X)^F\Def\Ker(H^1(X,\scrO_X)\xrightarrow{F-{\rm id}} H^1(X,\scrO_X))$. Then the Serre duality $H^0(X,\Omega^1_X)\simeq H^1(X,\scrO_X)^{\vee}$ induces an isomorphism of $\F_p$-vector spaces
\begin{equation}\label{eq:duality X}
H^1_{\fppf}(X,\mu_p)\simeq H^0(X,\Omega^1_X)^{\mathcal{C}}\simeq (H^{1}(X,\scrO_X)^{F})^{\vee}\simeq H^1_{\et}(X,\Z/p\Z)^{\vee}.
\end{equation}

\begin{lem}(cf.\ \cite[Proposition 2.5]{Pa95}\cite[Lemma 2.3]{Bo01})\label{lem:fq et fg}
Let $Y\to X$ be a connected finite \'etale Galois covering over $X$ with Galois group $H\Def\Gal(Y/X)$ having prime-to-$p$ order. Then there exists a bijection between the set of isomorphism classes of connected finite \'etale Galois covering $Z\to X$ which dominates the covering $Y\to X$ whose Galois group $\Gal(Z/X)$ is isomorphic to an extension of $H$ by an elementary abelian $p$-group and the set of $H$-submodules of $\Hom(\pi_1^{\et}(Y),\Z/p\Z)$. 
\end{lem}

Note that the isomorphism (\ref{eq:duality X}) for the curve $Y$ is compatible with the Galois module structure, hence we have an isomorphism of $H$-modules~(cf.\ \cite[\S2]{Bo01}),
\begin{equation}\label{eq:duality Y}
J_Y[p](k)\simeq\Hom(\pi^{\et}_1(Y),\Z/p\Z)^{\vee}.\end{equation}

We now consider the problem describing the set of finite quotients $\pi^{\et}_A(X)^{\scrC}$~(cf.\ (\ref{eq:Pi^C})). 
Recall that the maximal pro-prime-to-$p$ quotient $\pi^{\et}_1(X)^{(p')}$ is isomorphic to the profinite group $\Pi^{(p')}_{g,0}$~(cf.\ Corollary \ref{cor:pi^t fin gen}). Therefore, Lemma \ref{lem:fq et fg} suggests that it remains to determine the $\Gal(Y/X)$-module structure of $\Hom(\pi^{\et}_1(Y),\Z/p\Z)$ for every prime-to-$p$ Galois covering $Y\to X$. In fact, the Galois module structure on $\Hom(\pi^{\et}_1(Y),\Z/p\Z)$ can be described in terms of the \textit{generalized Hasse--Witt invariants}. For the detail, see \cite{Pa95}.


\section{The tame fundamental group scheme}\label{sec:tame fgs}

\subsection{Finite linearly reductive group schemes}\label{sec:lin red}

Let $S$ be a scheme. For an affine flat $S$-group scheme $G$, we will denote by $\X(G)$ the group of \textit{characters} of $G$, i.e.\ $\X(G)\Def\Hom_{S\text{-gr}}(G,\G_{m,S})$. 

For any  abelian group $A$, we denote by $\Diag_S(A)$ the \textit{diagonalizable} $S$-group scheme associated with $A$~\cite[Section 2.2]{Wa79}. For example, we have $\Diag_S(\Z)=\G_{m,S}$ and $\Diag_S(\Z/m\Z)=\mu_{m,S}~(m\in\Z)$. 
An affine flat $S$-group scheme $G$ is said to be \textit{diagonalizable} if it is isomorphic to the diagonalizable $S$-group scheme $\Diag_S(A)$ for some abelian group $A$. Then the correspondence $A\mapsto\Diag_S(A)$ gives an anti-equivalence of categories between the category of abelian groups and the category of diagonalizable $S$-group schemes. A quasi-inverse functor is given by taking the groups of characters, $G\mapsto\X(G)$. Furthermore, this equivalence of categories is compatible with any base change. Namely, for any morphism $T\to S$ and any abelian group $A$, we have $\Diag_S(A)\times_S T\simeq\Diag_T(A)$. Note that if $A$ is finite, then the Cartier dual $G^D=\uHom_S(G,\G_{m,S})$ of the diagonalizable $S$-group scheme $G=\Diag_S(A)$ is canonically isomorphic to the constant $S$-group scheme $\underline{A}_S$ associated with the finite abelian group $A$.

A finite flat $S$-group scheme $G$ is said to be \textit{linearly reductive} if the functor $\Qcoh(\cB_S G)\to\Qcoh(S);~\calF\mapsto\calF^G$ 
is exact~(cf.~\cite[Definition 2.2]{AOV08}), where $\calF^G$ denotes the $G$-invariant subsheaf of $\calF$. If $S=\Spec k$ is the spectrum of a field $k$, then the condition can be replaced by the condition that the category $\Rep(G)$ is semi-simple. The class of linearly reductive group schemes is stable under any base change $S'\to S$ and admits faithfully flat descent~(cf.\ \cite[Proposition 2.4]{AOV08}). Moreover, the class of linearly reductive group schemes is closed under taking subgroup schemes, quotients and extensions~(cf.\ \cite[Proposition 2.5]{AOV08}). 

Now let us recall the following classification result. 

\begin{prop}(cf.\ \cite[Lemma 2.17]{AOV08})\label{prop:class lr}
Let $S$ be the spectrum of a strictly henselian local ring and $G$ a finite linearly reductive $S$-group scheme. Then there exists a diagonalizable normal subgroup scheme $\Delta$ of $G$ such that $G/\Delta$ is a constant $S$-group scheme of order invertible in $S$.
\end{prop}

Let $S=\Spec k$ be the spectrum of a field $k$. If $k$ is of characteristic $0$, then any finite $k$-group scheme $G$ is linearly reductive. So let us assume that $k$ is of positive characteristic $p>0$. If $k$ is algebraically closed, then the classification is quite simple. Let $G$ be a finite linearly reductive group scheme over an algebraically closed field $k$ of characteristic $p>0$. Then the connected-to-\'etale exact sequence~(cf.\ \cite[\S6.7]{Wa79})
\begin{equation*}
1\to G^0\to G\to \pi_0(G)\to 1
\end{equation*}
admits a unique section $s\colon\pi_0(G)\to G$, which is in fact induced by the reduced closed subgroup scheme $G_{\rm red}$ of $G$, i.e.\ the composition $G_{\rm red}\hookrightarrow G\to\pi_0(G)$ is an isomorphism of finite \'etale group schemes. As $\pi_0(G)$ is \'etale and linearly reductive, it must be isomorphic to a constant $k$-group scheme $\underline{H}$ associated with a finite group $H$ of prime-to-$p$ order. On the other hand, $G^0$ is connected and linearly reductive, it is isomorphic to the diagonalizable $k$ group scheme $\Diag(A)$ for some abelian $p$-group $A$. Therefore, $G$ is isomorphic to a semi-direct product $\underline{H}\ltimes\Diag(A)$. Thus, the isomorphism class of $G$ is uniquely determined by the groups $H, A$ and the conjugacy action $\underline{H}\to\uAut(\Diag(A))=\uAut(\underline{A})$, where the last identification is due to the fact that $\Diag(A)^D=\underline{A}$. However, as $k$ is algebraically closed, the \'etale $k$-group scheme $\uAut(\underline{A})$ is constant and the homomorphism $\underline{H}\to\uAut(\underline{A})$  is uniquely determined by the induced homomorphism between the groups of $k$-valued points $H\to\Aut(A)$. In particular, any finite linearly reductive $k$-group scheme $G$ is defined over an algebraic closure $\overline{\F}_p$ of the prime field $\F_p$. 

As a consequence, we have the following.  

\begin{prop}\label{prop:lr alg closed bc}
Let $K/k$ be an extension of algebraically closed fields of characteristic $p>0$. Then there exists an equivalence of categories between the category of finite linearly reductive $k$-group schemes and the category of finite linearly reductive $K$-group schemes. 
\end{prop}

Now we introduce a subclass of finite linearly reductive group schemes. 

\begin{definition}\label{def:category D}
Let $k$ be an algebraically closed field of characteristic $p>0$. A finite linearly reductive $k$-group scheme $G$ is said to be \textit{elementary} if the character group $\X(G^0)$ of the connected part $G^0$ is an elementary abelian $p$-group, or equivalently if $G^0$ is isomorphic to a finite direct product of $\mu_p$, i.e.\ $G^0\simeq\mu_{p,k}^{s}$ for some integer $s\ge 0$. We denote by $\scrD$ the category of finite elementary linearly reductive $k$-group schemes.  
\end{definition}

Then the above discussion implies the following result.

\begin{prop}\label{prop:cat C=D}
There exists a canonical equivalence of categories between the category $\scrC$ in Definition \ref{def:category C} and the category $\scrD$ of elementary linearly reductive $k$-group schemes. 
\end{prop}

\begin{proof}
For any $\Gamma\in\scrC$, we define the group scheme $\Psi(\Gamma)\in\scrD$ as follows. By definition, there exists a unique elementary abelian normal $p$-Sylow subgroup $P\triangleleft G$ so that $G\xrightarrow{\simeq}(G/P)\ltimes P$. Then we define $\Psi(\Gamma)\Def\underline{(G/P)}\ltimes\Delta(P^{\vee})$, which clearly belongs to $\scrD$. Then the correspondence $\scrC\to\scrD;G\mapsto \Psi(G)$ gives a desired equivalence of categories $\Psi\colon\scrC\xrightarrow{\simeq}\scrD$. 
\end{proof}

Finally let us recall the following two results due to Olsson.

\begin{prop}(cf.\ \cite[Corollary 4.3]{Ol12})\label{prop:Olsson 1}
Let $S$ be the spectrum of a discrete valuation ring with generic point $\eta\in S$, $X$ a smooth $S$-scheme, and $G$ a finite linearly reductive $S$-group scheme. Then for any $G_{\eta}$-torsor  $P_{\eta}\to X_{\eta}$, after taking a finite extension $S'\to S$ there exists an extension of $P_{\eta}$ to a $G$-torsor $P\to X$ and such an extension $P\to X$ is unique up to isomorphism. 
\end{prop}

\begin{prop}(cf.\ \cite[Proposition 4.5]{Ol12})\label{prop:Olsson 2}
Let $k$ be a separably closed field. Let $X$ be a smooth scheme over $k$ and $\bfD=(D_i)_{i=1}^n$ a family of distinct reduced irreducible effective Cartier divisors $D_i$ on $X$, and $U\Def X\setminus D$, where $D=\cup_{i=1}^n\Supp(D_i)$. Let $G$ be a finite linearly reductive $k$-group scheme. Then for any $G$-torsor $P\to U$, there exists an $n$-tuple  $\bfr=(r_i)_{i=1}^n$ of positive integers $r_i$ such that $P\to U$ extends to a $G$-torsor $\overline{P}\to\sqrt[\bfr]{\bfD/X}$ over the root stack associated with the data $(X,\bfD,\bfr)$~(cf.\ \ref{para:def of root stack}). Moreover, the  extension over $\sqrt[\bfr]{\bfD/X}$ is unique up to unique isomorphism.   
\end{prop}

Note that if $G$ is constant of order prime to the characteristic of $k$, then the composite map $\overline{P}\to \sqrt[\bfr]{\bfD/X}\to X$ is nothing other than the normalization of the \'etale covering $P\to U$.


\subsection{The tame fundamental group scheme}\label{sec:tame fgs def}

In this subsection, we recall the definition of the \textit{tame fundamental group scheme}, which is the main object in the present paper.

\begin{definition}(cf.~\cite[Definition 10.1]{BV15})
A finite stack $\Gamma$ over a field $k$~(cf.\ \ref{para:fin stack}) is said to be \textit{tame} if the global section functor $\Coh(\Gamma)\to\Vecf_k;F\mapsto H^0(\Gamma,F)$ is exact.
\end{definition}

For example, if $G$ is a finite $k$-group scheme, then the classifying stack $\cB G$ is tame if and only if $G$ is linearly reductive~(cf.\ \S\ref{sec:lin red}). 

\begin{prop}(cf.~\cite[\S10]{BV15})\label{prop:def tame pi}
Let $\calX$ be an inflexible algebraic stack of finite type over a field $k$~(cf.\ \ref{para:Nori FG}). Then there exists a profinite tame gerbe $\Pi$ together with a morphism $\calX\to\Pi$ such that for any finite tame stack $\Gamma$ over $k$, the induced functor
\begin{equation*}
\Hom_k(\Pi,\Gamma)\lto\Hom_k(\calX,\Gamma)
\end{equation*}
is an equivalence of categories.
\end{prop}

\begin{definition}(cf.~\cite[Definition 10.4]{BV15})
With the above notation, the profinite tame gerbe $\Pi$ in Proposition \ref{prop:def tame pi} is unique up to unique isomorphism for $\calX/k$. We denote it by $\Pi^{\tame}_{\calX/k}$ and call it the \textit{tame fundamental gerbe} for $\calX$ over $k$. Moreover, if $\calX$ admits a $k$-rational point $x:\Spec k\to\calX$, then we denote by $\pi^{\tame}(\calX,x)$ the automorphism group scheme $\underline{\mathrm{Aut}}_k(\xi)$ of the object $\xi:\Spec k\xrightarrow{x}\calX\to\Pi^{\tame}_{\calX/k}$ and call it the \textit{tame fundamental group scheme} for $(\calX,x)$.
\end{definition}

\begin{rem}\label{rem:N vs tame}
Let $\pi^{\N}(\calX,x)$ be the fundamental group scheme for $(\calX,x)$~(cf.\ \ref{para:Nori fgs}, see also \cite[Chapter II]{No82}). Then the tame fundamental group scheme $\pi^{\tame}(\calX,x)$ is canonically isomorphic to the maximal pro-linearly reductive quotient of $\pi^{\N}(\calX,x)$. 
\end{rem}

As $\Pi^{\tame}_{\calX/k}$ is a profinite gerbe, the morphism $\calX\to\Pi^{\tame}_{\calX/k}$ factors through the Nori fundamental gerbe $\Pi^{\N}_{\calX/k}$~(cf.\ \ref{para:Nori FG}), 
\begin{equation*}
\begin{xy}
\xymatrix{
\calX\ar[r]\ar[rd]&\Pi^{\N}_{\calX/k}\ar[d]\\
&\Pi^{\tame}_{\calX/k}
}
\end{xy}
\end{equation*}
and the resulting morphism $\Pi^{\N}_{\calX/k}\to\Pi^{\tame}_{\calX/k}$ gives a gerbe. 
Recall that if $\calX$ is proper over $k$, then we have an equivalence of tannakian categories over $k$,
\begin{equation*}
\Vect(\Pi^{\N}_{\calX/k})\xrightarrow{\simeq}\EFin(\calX),
\end{equation*}
where $\EFin(\calX)$ is the full subcategory of $\Vect(\calX)$ consisting of essentially finite bundles~(cf.\ \ref{para:EFin}). An essentially finite bundle $E$ on $\calX$ is said to be \textit{tamely finite} if all the indecomposable components of all the tensor powers $E^{\otimes n}$ are irreducible~(cf.~\cite[Definition 12.1]{BV15}). We now define the category $\TFin(\calX)$ to be the full tannakian subcategory of $\EFin(\calX)$ which consists of tamely finite bundles on $\calX$. Then, the fully faithful tensor functor 
\begin{equation*}
\Vect(\Pi^{\tame}_{\calX/k})\hookrightarrow\Vect(\Pi^{\N}_{\calX/k})\xrightarrow{~\simeq~}\EFin(\calX)\subset\Vect(\calX)
\end{equation*}
induces an equivalence of categories (cf.\,\cite[Theorem 12.2]{BV15})
\begin{equation}\label{eq:TFin}
\Vect(\Pi^{\tame}_{\calX/k})\xrightarrow{~\simeq~}\TFin(\calX).
\end{equation}
In particular, $\TFin(\calX)$ is a tannakian category over $k$ and it is the largest tannakian semi-simple subcategory of $\EFin(\calX)$.

\subsection{Structural results}\label{sec:str tame fg}

In this subsection, we will see several structural results for the tame fundamental group schemes of curves. All the results should be well-understood by the experts, but we put them here for lack of references.

Let us begin with the following version of \cite[Theorem I]{ABETZ}.

\begin{prop}\label{prop:str tame pi proper}
Let $\calX$ be a proper inflexible algebraic stack over a perfect field $k$ of characteristic $p>0$. Let $f:\calY\to\calX$ be a Nori-reduced $G$-torsor~(cf.\ \ref{para:Nori-reduced}), where $G$ is a finite \'etale linearly reductive $k$-group scheme. Then all the squares in the following diagram are 2-Cartesian.
\begin{equation*}
\begin{xy}
\xymatrix{
\calY\ar[r]\ar[d]_{f}&\Pi^{\tame}_{\calY/k}\ar[r]\ar[d]_{f_*}&\Spec k\ar[d]\\
\calX\ar[r]&\Pi^{\tame}_{\calX/k}\ar[r]&\cB G.
}
\end{xy}
\end{equation*} 
\end{prop}

\begin{proof}
We shall adapt the argument in the proof of \cite[Theorem I]{ABETZ}.
We define a finite stack $\Pi$ over $\Pi^{\tame}_{\calX/k}$ to be the 2-fiber product
\begin{equation*}
\begin{xy}
\xymatrix{\ar@{}[rd]|{\square}
\Pi\ar[r]\ar[d]_{f_{\star}}&\Spec k\ar[d]\\
\Pi^{\tame}_{\calX/k}\ar[r]&\cB G.
}
\end{xy}
\end{equation*}
Hence, all the squares in the diagram 
\begin{equation*}
\begin{xy}
\xymatrix{
\calY\ar[r]^{v~~}\ar[d]_{f}&\Pi\ar[r]\ar[d]_{f_{\star}}&\Spec k\ar[d]\\
\calX\ar[r]_{u~~}&\Pi^{\tame}_{\calX/k}\ar[r]&\cB G.
}
\end{xy}
\end{equation*} 
are 2-Cartesian. It suffices to show that the morphism $v\colon\calY\to\Pi$ gives the tame fundamental gerbe of $\calY$. First note that since $f:\calY\to\calX$ is Nori-reduced, the morphism $\Pi^{\tame}_{\calX/k}\to\cB G$ is a gerbe and so is $\Pi\to\Spec k$. 
Moreover, since $\Pi\to\Pi^{\tame}_{\calX/k}$ is representable and $\Pi^{\tame}_{\calX/k}$ is tame, the gerbe $\Pi$ is tame as well. In particular, by \cite[Proposition 10.3]{BV15}, the category $\Vect(\Pi)$ is semisimple.

Since $\calX$ is proper over $k$, the pullback functor $u^*$ induces an equivalence of tannakian categories (cf.\,(\ref{eq:TFin}))
\begin{equation*}
u^*:\Vect(\Pi^{\tame}_{\calX/k})\xrightarrow{~\simeq~}\TFin(\calX),
\end{equation*} 
hence $u_*\scrO_{\calX}\simeq\scrO_{\Pi^{\tame}_{\calX/k}}$ (cf.\,\cite[Lemma 1.22]{ABETZ}). As $f_{\star}$ is faithfully flat, by flat base change theorem, we also have $v_*\scrO_{\calY}\simeq\scrO_{\Pi}$, which implies that the pullback functor
\begin{equation}\label{eq:v^*}
v^*:\Vect(\Pi)\to\Vect(\calY)
\end{equation}  
is fully faithful. However, as $\Vect(\Pi)$ is a semisimple tannakian category, the essential image of $v^*$ must be contained in $\TFin(\calY)$. To prove that $\calY\to\Pi$ gives the tame fundamental gerbe of $\calY$, it suffices to show that the functor (\ref{eq:v^*}) has essential image $\TFin(\calY)$.

By applying \cite[Lemma 2.7]{ABETZ} to the diagram
\begin{equation*}
\begin{xy}
\xymatrix{
\Vect(\calY)&\TFin(\calY)\ar@{_{(}->}[l]&\Vect(\Pi)\ar@{_{(}->}[l]\ar@/_8mm/[ll]_{v^*}\\
\Vect(\calX)\ar[u]^{f^*}&\TFin(\calX)\ar@{_{(}->}[l]\ar[u]&\Vect(\Pi^{\tame}_{\calX/k})\ar[l]_{\simeq}\ar[u]\ar@/^8mm/[ll]_{u^*},
}
\end{xy}
\end{equation*}
we can conclude that for any $V\in\Vect(\calY)$, it is contained in the essential image of $v^*$ if and only if $f_*V$ is tamely finite. Therefore, to complete the proof, we have only 
to show that for any tamely finite bundle $V\in\TFin(\calY)$, the pushforward sheaf $f_*V$ is tamely finite on $\calX$.

This is a consequence of the existence of Galois envelopes in the sense of (cf.\,\cite[Definition 3.8]{ABETZ}). In our situation, this ensures that for any $H$-torsor $\calZ\to\calY$ with $H$ linearly reductive, there exists a $\Gamma$-torsor $\calP\to\calX$ where $\Gamma$ is linearly reductive together with homomorphisms $\alpha:\Gamma\to G$ and $\beta:\Ker(\alpha)\to H$ and a commutative diagram
\begin{equation*}
\begin{xy}
\xymatrix{
\calP\ar[d]\ar[dr]_{g}\ar[drr]^{\pi}&&\\
\calZ\ar[r]_h&\calY\ar[r]_{f}&\calX,
}
\end{xy}
\end{equation*}
where $g\colon \calP\to\calY$ is $\Ker(\alpha)$-equivariant and $\calP\to\calZ$ is $\Ker(\beta)$-equivariant.

For any tamely finite bundle $V\in\TFin(\calY)$, by definition, there exists an $H$-torsor $h:\calZ\to\calY$, where $H$ is a finite linearly reductive $k$-group scheme, such that $h^*V\simeq\scrO_{\calZ}^{\oplus m}$. Therefore, if $\pi:\calP\to\calX$ denotes the Galois envelope associated with the tower $\calZ\to\calY\to\calX$  as above and $g:\calP\to\calY$ the $\Ker(\alpha)$-equivariant morphism, then $V\subseteq g_*g^*V\simeq g_*\scrO_{\calP}^{\oplus m}$, hence $f_*V\subseteq f_*g_*\scrO_{\calP}^{\oplus m}\simeq\pi_*\scrO_{\calP}^{\oplus m}$, where the latter vector bundle is tamely finite because $\pi:\calP\to\calX$ is a finite linearly reductive torsor. This completes the proof.
\end{proof}

Now we study the tame fundamental group schemes of curves. Let $k$ be an algebraically closed field $k$ of characteristic $p>0$. 
Let $X$ be a projective smooth curve over $k$. Let $S$ be a finite (possibly empty) set of closed points of $X$ with $n\Def\#S$. Let $U\Def X\setminus S$. We set $\bfD\Def (x)_{x\in S}$. For any $n$-tuple $\bfr=(r_x)_{x\in S}$ of positive integers, we define
\begin{equation}\label{eq:fX}
\fX^{\bfr}\Def\sqrt[\bfr]{\bfD/X}
\end{equation}
to be the root stack associated with the data $(X,\bfD,\bfr)$ (cf.\,\ref{para:def of root stack}).

\begin{prop}\label{prop:cor of Olsson 2}
With the above notation, there exists a natural isomorphism of affine gerbes over $k$,
\begin{equation*}
\Pi^{\tame}_{U/k}\xrightarrow{~\simeq~}\varprojlim_{\bfr}\Pi^{\tame}_{{\fX}^{\bfr}/k}
\end{equation*} 
\end{prop}

\begin{proof}
The morphisms into root stacks $U\to\fX^{\bfr}$ induce morphisms between the tame fundamental gerbes $\Pi^{\tame}_{U/k}\to\Pi^{\tame}_{\fX^{\bfr}/k}$, hence we obtain a morphism $\Pi^{\tame}_{U/k}\to\varprojlim_{\bfr}\Pi^{\tame}_{\fX^{\bfr}/k}$. This is surjective. The injectivity follows from Proposition \ref{prop:Olsson 2}.
\end{proof}

\begin{cor}\label{cor:Olsson 2 open vs proper}
With the above notation, we further fix a $k$-rational point $x_0\in U$. Then for each $x\in S$, there exists a homomorphism $\delta_x\colon\Diag(\Q/\Z)\to\pi^{\tame}(U,x_0)$ which is canonical up to conjugation so that the kernel of the surjective homomorphism $\pi^{\tame}(U,x_0)\twoheadrightarrow\pi^{\tame}(X,x_0)$ is the normal subgroup generated by the images of $\delta_x$. 
\end{cor}

\begin{proof}
We will use the same notation as in Proposition \ref{prop:cor of Olsson 2}. For each $\bfr=(r_x)_{x\in S}$, let $\coprod_{x\in S}\cB\mu_{r_x}\to\fX^{\bfr}$ be the coproduct of the maps from the residual gerbes into the root stack $\fX^{\bfr}$~(cf.\ \ref{para:res gerbe}). By taking the limit, we get a sequence of morphisms of pro-algebraic stacks
\begin{equation*}
\coprod_{x\in S}\cB\Delta(\Q/\Z)\to\varprojlim_{\bfr}\fX^{\bfr}\to\varprojlim_{\bfr}\Pi_{\fX^{\bfr}/k}^{\tame}\to\Pi^{\tame}_{X/k}.
\end{equation*}
If $\xi\colon\Spec k\to\Pi^{\tame}_{X/k}$ is any section, then by Biswas--Borne's theorem~(cf.\ \cite[Corollary 3.6]{BB17}), the map $\varprojlim_{\bfr}\Pi^{\tame}_{\fX^{\bfr}/k}\to\Pi^{\tame}_{X/k}$ is a relative gerbe whose fiber 
\begin{equation*}
(\varprojlim\Pi^{\tame}_{\fX^{\bfr}/k})\times_{\Pi^{\tame}_{X/k}}{{}_{\xi}\Spec k}
\end{equation*}
over $\xi$ is generated by the image of the map $\coprod_{x\in S}{\cB}\Delta(\Q/\Z)\to\varprojlim_{\bfr}\Pi^{\tame}_{\fX^{\bfr}/k}$. The assertion is then a consequence of Proposition \ref{prop:cor of Olsson 2}.
\end{proof}

\begin{cor}\label{cor:Olsson 2 mu_p}
With the above notation, suppose that $n>0$ and fix a point $x_0\in S$. Then for any integer $m>0$, there exists an exact sequence of $\Z/p^m\Z$-modules,
\begin{equation*}
0\to H^1_{\fppf}(X,\mu_{p^m})\to H^{1}_{\fppf}(U,\mu_{p^m})\to \oplus_{x\in S\setminus\{x_0\}}\Z/p^m\Z\to 0.
\end{equation*} 
\end{cor}

\begin{proof}
First consider the case $n=1$. In this case, the natural restriction map $\Pic(X)\twoheadrightarrow\Pic(U)$ induces an isomorphism $\Pic(X)[p^m]\xrightarrow{\simeq}\Pic(U)[p^m]$. Thus, 
by the Kummer theory, we obtain a natural isomorphism $H^1_{\fppf}(X,\mu_{p^m})\xrightarrow{\simeq}H^1_{\fppf}(U,\mu_{p^m})$. Therefore, the assertion is true for $n=1$. 

Next we will deal with the arbitrary case. Let us  consider the one-punctured curve $X_0\Def X\setminus\{x_0\}$. By the previous case $n=1$, we are reduced to show that there exists an exact sequence of $\Z/p^{m}\Z$-modules,
\begin{equation*}
0\to H^1_{\fppf}(X_0,\mu_{p^m})\to H^{1}_{\fppf}(U,\mu_{p^m})\to \oplus_{x\in S\setminus\{x_0\}}\Z/p^m\Z\to 0.
\end{equation*} 
Let $\mathbf{p}^m\Def (p^m)_{x\in S\setminus\{x_0\}}$. By Proposition \ref{prop:Olsson 2}, we have an isomorphism of $\Z/p^m\Z$-modules,
\begin{equation*}
H^1_{\fppf}(\fX_0^{\mathbf{p}^m},\mu_{p^m})\xrightarrow{\simeq}H^1_{\fppf}(U,\mu_{p^m}),
\end{equation*}
where $\fX_0^{\mathbf{p}^{m}}$ is the root stack associated with the data $(X_0,S\setminus\{x_0\},\mathbf{p}^{m})$.
Therefore, it suffices to show that
there exists an exact sequence of $\Z/p^m\Z$-modules,
\begin{equation}\label{eq:Olsson 2 mu_p}
0\to H^1_{\fppf}(X_0,\mu_{p^m})\to H^{1}_{\fppf}(\fX_0^{\mathbf{p}^m},\mu_{p^m})\to \oplus_{x\in S\setminus\{x_0\}}\Z/p^m\Z\to 0. 
\end{equation} 
Recall that there exists an exact sequence of abelian groups
\begin{equation*}
0\to\Pic(X_0)\to\Pic(\fX_0^{\mathbf{p}^m})\to\oplus_{x\in S\setminus\{x_0\}}\Z/p^m\Z\to 0
\end{equation*}
(cf.\ \cite[\S5.4]{Br09}). Therefore, by the Kummer theory together with the snake lemma, we obtain the desired exact sequence (\ref{eq:Olsson 2 mu_p}). This completes the proof.  
\end{proof}

Finally we notice the following fact.

\begin{prop}\label{prop:str tame ger}
With the same notation as above, let $f:V\lto U$ be a Nori-reduced $G$-torsor (cf.\,\ref{para:Nori-reduced}), where $G$ is a constant $k$-group scheme of prime-to-$p$ order. Then all the squares in the following diagram are 2-Cartesian.
\begin{equation*}
\begin{xy}
\xymatrix{
V\ar[r]\ar[d]_{f}&\Pi^{\tame}_{V/k}\ar[r]\ar[d]_{f_*}&\Spec k\ar[d]\\
U\ar[r]&\Pi^{\tame}_{U/k}\ar[r]&\cB G.
}
\end{xy}
\end{equation*} 
\end{prop}

\begin{proof}
The assertion is immediate from Propositions  \ref{prop:str tame pi proper} and \ref{prop:cor of Olsson 2}. 
\end{proof}

As a consequence, we obtain the following version of \cite[Theorem 2.9]{EHS08}.

\begin{cor}\label{cor:tame fgs fet}
Let $k$ be an algebraically closed field of characteristic $p>0$. Let $U$ be a smooth connected curve over $k$ with $x\in U(k)$ a $k$-rational point. Let $H$ be a finite group of prime-to-$p$ order. Suppose given a surjective homomorphism $\phi\colon\pi^{\tame}(U,x)\twoheadrightarrow\underline{H}$ onto the constant group scheme $\underline{H}$. Let $(V,y)\to(U,x)$ be the corresponding pointed $\underline{H}$-torsor. Then we have an exact sequence of affine $k$-group schemes
\begin{equation*}
1\to\pi^{\tame}(V,y)\to\pi^{\tame}(U,x)\xrightarrow{\phi}\underline{H}\to 1.
\end{equation*}
\end{cor}

\begin{proof}
This directly follows from Proposition \ref{prop:str tame ger}.
\end{proof}

\subsection{The base change theorem}\label{sec:bc tame fgs}

Let $X$ be a proper connected and reduced scheme over an algebraically closed field $k$ and $K$ an algebraically closed extension of $k$. Let $x\in X(k)$ be a $k$-rational point. By \cite[Expos\'e X, Corollaire 1.8]{Gr71}, the natural homomorphism $\pi^{\et}_1(X_K,x_K)\to\pi^{\et}_1(X,x)$ of the \'etale fundamental groups is an isomorphism. However, it is known that the natural map between Nori's fundamental group schemes
\begin{equation}\label{eq:h^N}
\pi^{\N}(X_K,x_K)\to\pi^{\N}(X,x)\times_k K
\end{equation}
is not an isomorphism in general (cf.\,\cite{MS}\cite{Pauly}). In this subsection, we will prove that the  homomorphism (\ref{eq:h^N}) induces an isomorphism of the tame fundamental group schemes for smooth curves.

Let $k$ be an algebraically closed field of characteristic $p>0$. Let $X$ be a proper smooth curve over $k$ and $S\subset X$ a finite (possibly empty) set of closed points. Let $U\Def X\setminus S$ and fix a $k$-rational point $x\in U(k)$.

\begin{prop}\label{prop:bc tame fgs}
With the above notation, let $K/k$ be an extension of algebraically closed fields. Then we have a natural isomorphism of affine $K$-group schemes
\begin{equation*}
h_K\colon \pi^{\tame}(U_K,x_K)\xrightarrow{\simeq}\pi^{\tame}(U,x)\times_k K.
\end{equation*}
\end{prop}

\begin{proof}
First let us show the surjectivity of the natural homomorphism $h_K\colon \pi^{\tame}(U_K,x_K)\to\pi^{\tame}(U,x)\times_k K$. Let $\phi\colon\pi^{\tame}(U,x)\twoheadrightarrow G$ be any finite quotient map. We have to show that the composition of the homomorphism
\begin{equation*}
\pi^{\tame}(U_K,x_K)\xrightarrow{h_K}\pi^{\tame}(U,x)\times_k K\xrightarrow{\phi\times_k K}G\times_k K
\end{equation*}
is surjective. Let $(P,p)\to (U,x)$ be the pointed $G$-torsor associated with $\phi$. We set $(V,y)\Def (P,p)/G^0$, which is a pointed $\pi_0(G)$-torsor. Thanks to Corollary \ref{cor:tame fgs fet}, we have the commutative diagram
\begin{equation*}
\begin{xy}
\xymatrix{
1\ar[d]&1\ar[d]&1\ar[d]\\
\pi^{\tame}(V_K,y_K)\ar[r]\ar[d]&\pi^{\tame}(V,y)\times_kK\ar@{->>}[r]\ar[d]&G^0\times_k K\ar[d]\\
\pi^{\tame}(U_K,x_K)\ar[r]^{h_K\quad}\ar[d]&\pi^{\tame}(U,x)\times_kK\ar@{->>}[r]^{\qquad\phi\times_kK}\ar[d]&G\times_k K\ar[d]\\
\pi_0(G)\times_kK\ar@{=}[r]\ar[d]&\pi_0(G)\times_kK\ar@{=}[r]\ar[d]&\pi_0(G)\times_k K,\ar[d]\\
1&1&1
}
\end{xy}
\end{equation*}
where all the vertical sequences are exact. Therefore, the map $(\phi\times_k K)\circ h_K$ is surjective if and only if so is the map
\begin{equation*}
\pi^{\tame}(V_K,y_K)\to\pi^{\tame}(V,y)\times_kK\twoheadrightarrow G^0\times_k K.
\end{equation*}
The latter condition is valid because we have an isomorphism
\begin{equation}\label{eq:bc diag}
H^1_{\fppf}(V,G^0)\xrightarrow{\simeq}H^1_{\fppf}(V_K,G^0\times_kK),
\end{equation}
which follows from Corollary \ref{cor:Olsson 2 mu_p} together with the isomorphism $\Pic_{X_K}^0\simeq\Pic^0_X\times_k K$. This completes the proof of the surjectivity of the map $h_K$. 

For the injectivity of the map $h_K$, it suffices to notice that any finite quotient map 
\[
\psi\colon \pi^{\tame}(U_K,x_K)\twoheadrightarrow G'
\] 
factors through the map $h_K$. Indeed, by Proposition \ref{prop:class lr}, we may assume that $G'=G\times_k K$ for some finite linearly reductive $k$-group scheme $G$. It suffices to show that the pointed $G\times_k K$-torsor $(P,p)\to (U_K,x_K)$ corresponding to $\psi$ is defined over $(U,x)$. This is well known when $G$ is \'etale, in which the torsor $P\to U_K$ is a prime-to-$p$ Galois covering, hence we are reduced to the connected case $G=G^0$, for which the claim immediately follows from the isomorphism (\ref{eq:bc diag}). This completes the proof of the proposition.    
\end{proof}


\section{The cospecialization map for the tame fundamental group scheme}\label{sec:main}

\subsection{Finite quotients of the tame fundamental group scheme}\label{sec:fq tame fgs}

Let $k$ be an algebraically closed field of characteristic $p>0$ and $X$ a proper smooth connected curve over $k$. Let $S\subset X(k)$ be a finite (possibly empty) set of closed points and set $U\Def X\setminus S$. Fix a $k$-rational point $x\in U(k)$. In this subsection, we follow the argument in \S\ref{sec:fq et fg} to investigate finite quotients of the tame fundamental group scheme $\pi^{\tame}(U,x)$.

Let $H$ be a fixed finite group of prime-to-$p$ order. 
Suppose given a surjective homomorphism 
\[
\phi\colon\pi^{\tame}(U,x)\twoheadrightarrow\underline{H}
\] 
and let $(V,y)\to (U,x)$ be the corresponding pointed $\underline{H}$-torsor.

\begin{lem}\label{lem:count dq}
There exists a bijection between the set of isomorphism classes of surjective homomorphisms $\pi^{\tame}(V,y)\twoheadrightarrow\mu_p^{s}$ for some integer $s\ge 0$ and the set of subspaces of the $\F_p$-vector space $\Hom(\pi^{\tame}(V,y),\mu_p)$.  
\end{lem}

Here, any two surjective homomorphisms $\psi_1\colon \pi^{\tame}(V,y)\twoheadrightarrow\mu_p^{s_1}$ and $\psi_2\colon \pi^{\tame}(V,y)\twoheadrightarrow\mu_p^{s_2}$ are said to be \textit{isomorphic} if $s_1=s_2$ and there exists an isomorphism of group schemes $\alpha\in\Aut(\mu_p^{s_1})$ such that $\psi_1=\alpha\circ \psi_2$. 

\begin{proof}
Let $S_1(\phi)$ be the set of isomorphism classes of surjective homomorphisms $\pi^{\tame}(V,y)\twoheadrightarrow\mu_p^s$ for some $s\ge 0$. Let $S_2(\phi)$ be the set of subspaces of the $\F_p$-vector space $\Hom(\pi^{\tame}(V,y),\mu_p)$. Then the map $M\colon S_1(\phi)\to S_2(\phi);\psi\mapsto M(\psi)$ defined by
\begin{equation}\label{eq:M psi}
M(\psi)\Def\im\bigl(\Hom(\mu_p^s,\mu_p)\xrightarrow{\psi^*}\Hom(\pi^{\tame}(V,y),\mu_p)\bigl)
\end{equation}
gives the desired bijective map.
\end{proof}

Note that the $\F_p$-vector space $\Hom(\pi^{\tame}(V,y),\mu_p)$ is canonically isomorphic to the cohomology group $H^1_{\fppf}(V,\mu_p)$ as an $\F_p[H]$-module. Let $M$ be a finitely generated $\F_p[H]$-module. 
Then the corresponding action $H\to\Aut(M)\Def\GL_M(\F_p)$ induces a homomorphism of $k$-group schemes $\underline{H}\to\uAut(\Diag(M))$. Hence, we obtain a finite linearly reductive $k$-group scheme $G\Def\underline{H}\ltimes\Diag(M)$.

\begin{lem}\label{lem:count lrq}
Suppose given a surjective homomorphism $\psi\colon\pi^{\tame}(V,y)\twoheadrightarrow\mu_p^s$ and let $M(\psi)\subset\Hom(\pi^{\tame}(V,y),\mu_p)$ be the corresponding subspace~(cf.\ Lemma \ref{lem:count dq} and (\ref{eq:M psi})). Let $\Gamma_{\psi}\vartriangleleft\pi^{\tame}(V,y)^{\rm ab}$ be the kernel of the homomorphism $\psi^{\ab}\colon\pi^{\tame}(V,y)^{\rm ab}\twoheadrightarrow\mu_p^s$ induced by $\psi$. Then the following conditions are equivalent. 
\begin{enumerate}
\renewcommand{\labelenumi}{(\alph{enumi})}
\item The subspace $M(\psi)$ is stable under the action by $H$ on $\Hom(\pi^{\tame}(V,y),\mu_p)$.

\item The subgroup scheme $\Gamma_{\psi}<\pi^{\tame}(V,y)^{\rm ab}$ is stable under the conjugacy action by $\underline{H}$.  

\item The surjective homomorphism $\psi\colon\pi^{\tame}(V,y)\twoheadrightarrow\mu_p^s$ extends to a surjective homomorphism $\widetilde{\psi}\colon\pi^{\tame}(U,x)\twoheadrightarrow G$ onto a semi-direct product $G\Def\underline{H}\ltimes\mu_p^s$ and the composition $\pi^{\tame}(U,x)\overset{\widetilde{\psi}}{\twoheadrightarrow} G\twoheadrightarrow G/\mu_p^s=\underline{H}$ coincides with $\phi\colon\pi^{\tame}(U,x)\twoheadrightarrow\underline{H}$. 
\end{enumerate}
\end{lem}

\begin{proof}
Thanks to Corollary \ref{cor:tame fgs fet}, we have an exact sequence
\begin{equation*}
1\to\pi^{\tame}(V,y)\to\pi^{\tame}(U,x)\to \underline{H}\to 1.
\end{equation*}
By the universality of the maximal abelian quotient $\pi^{\tame}(V,y)\twoheadrightarrow\pi^{\tame}(V,y)^{\ab}$, the kernel of the quotient map is normal in the group scheme $\pi^{\tame}(U,x)$. Hence, we get the quotient
\begin{equation*}
\pi^{\tame}(U,x)'\Def\pi^{\tame}(U,x)/\Ker(\pi^{\tame}(V,y)\twoheadrightarrow\pi^{\tame}(V,y)^{\ab}),
\end{equation*}
which fits into the commutative diagram with exact rows
\begin{equation*}
\begin{xy}
\xymatrix{
1\ar[r]&\pi^{\tame}(V,y)\ar[r]\ar@{->>}[d]&\pi^{\tame}(U,x)\ar[r]\ar[d]&\underline{H}\ar[r]\ar@{=}[d]&1\\
1\ar[r]&\pi^{\tame}(V,y)^{\ab}\ar[r]\ar@{->>}[d]&\pi^{\tame}(U,x)'\ar[r]&\underline{H}\ar[r]&1.\\
&\mu_p^s&&
}
\end{xy}
\end{equation*}
The condition (c) is then equivalent to the existence of the pushout $G$ as in the diagram
\begin{equation*}
\begin{xy}
\xymatrix{
1\ar[r]&\pi^{\tame}(V,y)^{\ab}\ar[r]\ar@{->>}[d]&\pi^{\tame}(U,x)'\ar[d]\ar[r]&\underline{H}\ar@{=}[d]\ar[r]&1.\\
1\ar[r]&\mu_p^s\ar[r]&G\ar[r]&\underline{H}\ar[r]&1
}
\end{xy}
\end{equation*}
Now the equivalences among the conditions (a), (b) and (c) clearly holds true. 
\end{proof}

The following is a version of Lemma \ref{lem:fq et fg}.

\begin{lem}\label{lem:fq tame fgs}
Let $H$ be a finite group of prime-to-$p$ order. Let $V\to U$ be a connected $\underline{H}$-torsor. Then the set of isomorphism classes of Nori-reduced $G$-torsors $P\to U$ (cf.\,\ref{para:Nori-reduced}) which dominate $V\to U$, where $G$ is an extension of $\underline{H}$ by a diagonalizable group scheme $\Delta$ with $p\X(\Delta)=0$ is in bijection with the set of $\F_p[H]$-submodules of $\Hom(\pi^{\tame}(V),\mu_p)$.
\end{lem}

\begin{proof}
This is a consequence of Lemmas \ref{lem:count dq} and \ref{lem:count lrq}. 
\end{proof}

\begin{prop}\label{prop:fq lr vs et}
Let $X$ be a proper smooth connected curve over $k$ and fix a $k$-rational point $x\in X(k)$. Let $H$ be a finite group of prime-to-$p$ order. Let $M$ be a finite dimensional $\F_p[H]$-module. Consider the finite group $\Gamma\Def H\ltimes (M^{\vee})$ and the finite linearly reductive $k$-group scheme $G\Def\underline{H}\ltimes\Delta(M)$. Then there exists a surjective homomorphism $\pi^{\tame}(X,x)\twoheadrightarrow G$ if and only if there exists a surjective homomorphism $\pi^{\et}_1(X,x)\twoheadrightarrow\Gamma$.  
\end{prop}

\begin{proof}
When $M=0$, the assertion is straightforward from the fact that $\pi^{\tame}(X,x)(k)=\pi^{\et}_1(X,x)^{(p')}$. Let us consider the arbitrary case. Suppose given a surjective homomorphism $\psi\colon\pi^{\tame}(X,x)\twoheadrightarrow G=\underline{H}\ltimes\Delta(M)$. Let $(Y,y)\to (X,x)$ be the Nori-reduced pointed $\underline{H}$-torsor associated with the surjective homomorphism $\pi^{\tame}(X,x)\overset{\psi}{\twoheadrightarrow}G\twoheadrightarrow\underline{H}$. Then, by Lemma \ref{lem:fq tame fgs}, the map $\psi$ makes $M$ a $H$-submodule of $\Hom(\pi^{\tame}(Y,y),\mu_p)=H^1_{\fppf}(Y,\mu_p)=J_Y[p](k)$. As there exists a natural isomorphism $J_Y[p](k)\simeq \Hom(\pi^{\et}_1(Y),\Z/p\Z)^{\vee}$ of $\F_p[H]$-modules~(cf.\ (\ref{eq:duality Y})), the dual $M^{\vee}$ can be embedded into $\Hom(\pi_1^{\et}(Y),\Z/p\Z)$ as an $H$-submodule. By Lemma \ref{lem:fq et fg}, this implies that $\Gamma$ appears as a quotient of $\pi^{\et}_1(X,x)$. The dual argument proves that the converse is also true. This completes the proof. 
\end{proof}

Finally we remark on the following finiteness result. 

\begin{prop}\label{prop:finiteness}
Let $U$ be a smooth connected curve over $k$ together with a $k$-rational point $x\in U(k)$. Then for any $G\in\scrD$~(cf.\ Definition \ref{def:category D}), the set $\Hom(\pi^{\tame}(U,x),G)$ of homomorphisms into $G$ is a finite set. 
\end{prop}

\begin{proof}
Let $H\Def G(k)=\pi_0(G)(k)$, which is a finite group of prime-to-$p$ order. As the set of subgroup schemes of $G$ is finite, it suffices to show that the set of surjective homomorphisms is a finite set. However, by Lemma \ref{lem:fq tame fgs}, this follows from the finiteness of the set of $H$-submodules of any finite dimensional $H$-module $M$. This completes the proof. 
\end{proof}


\subsection{Lifting problems for linearly reductive torsors}\label{sec:LP}

Let $k$ be an algebraically closed field of characteristic $p>0$. Let $R\Def k[[t]]$ be the ring of formal power series with coefficients in $k$ and $K\Def\Frac R=k((t))$ its field of fractions. Let us fix an algebraic closure $\overline{K}$ of $K$ and $\overline{R}$ the integral closrue of $R$ in $\overline{K}$. 
Let $S=\Spec R$.

Let $f\colon X\to S$ be a proper smooth $S$-curve with geometrically connected fibers. Let $X_{\olK}$ (respectively $X_k$) be the geometrically generic fiber (respectively the special fiber) of $X/S$. Let $D$ be a relatively \'etale Cartier divisor on $X/S$ of degree $n$. We set $U\Def X\setminus \Supp(D)$ and denote by $U_{\olK}$ (respectively by $U_k$) the geometrically generic fiber (respectively the special fiber) of the $S$-curve $U$. 

\begin{lem}\label{lem:LP^1 diag}
There exists a canonical surjective homomorphism between the cohomology groups $H^1_{\fppf}(U_{\overline{K}},\mu_{p})\twoheadrightarrow H^1_{\fppf}(U_k,\mu_{p})$.
\end{lem}

\begin{proof}
Thanks to Proposition \ref{prop:Olsson 1} applied to the diagonalizable group scheme $G=\mu_p$, we get a natural map
\begin{equation*}
H^1_{\fppf}(U_{\olK},\mu_p)\xleftarrow{\simeq}H^1_{\fppf}(U_{\olR},\mu_p)\to H^1_{\fppf}(U_k,\mu_p).
\end{equation*}
It suffices to prove that it is surjective. By Corollary \ref{cor:Olsson 2 mu_p}, one can reduce the problem to the proper case $X=U$. 

Let $\Pic^0_{X/S}$ be the connected component of the identity on the Picard scheme $\Pic_{X/S}$ for $X/S$. As $\Pic^0_{X/S}$ is an abelian $S$-scheme, we have the surjective reduction map
\begin{equation*}
\Pic^0(X_{\olK})=\Pic^0_{X/S}(\olK)=\Pic^0_{X/S}(\olR)\twoheadrightarrow\Pic^0_{X/S}(k)=\Pic^0(X_{k}).
\end{equation*}
As both the groups $\Pic^0(X_{\olK})$ and $\Pic^0(X_k)$ are divisible, the reduction map induces a surjective homomorphism
\begin{equation*}
\Pic^0(X_{\olK})[p]\twoheadrightarrow\Pic^0(X_k)[p].
\end{equation*}
On the other hand, by the Kummer theory, we have the canonical isomorphisms
\begin{equation*}
H^1_{\fppf}(X_{\olK},\mu_p)\xrightarrow{\simeq}\Pic^0(X_{\olK})[p]\quad\text{and}\quad H^1_{\fppf}(X_{k},\mu_p)\xrightarrow{\simeq}\Pic^0(X_{k})[p],
\end{equation*}
which implies that the map $H^1_{\fppf}(X_{\olK},\mu_p)\to H^1_{\fppf}(X_k,\mu_p)$ is surjective. This completes the proof. 
\end{proof}

\begin{prop}\label{prop:LP^1}
Let $G\in\scrD$~(cf.\ Definition \ref{def:category D}). Then any pointed $G$-torsor $(P,p)\to (U_k,x_k)$ can be lifted to a pointed $G_{\olK}$-torsor $(P_{\olK},p_{\olK})\to (U_{\olK},x_{\olK})$. \end{prop}

\begin{proof}
Let $H$ be a finite group of prime-to-$p$ order and $M$ a finitely generated $\F_p[H]$-module so that $G\simeq\underline{H}\ltimes\Diag(M)$. 
Without loss of generality, we may assume that $(P,p)\to (U_k,x_k)$ is a Nori-reduced $G$-torsor (cf.\,\ref{para:Nori-reduced}). Let $(V_k,y_k)\Def (P,p)/\Delta(M)$ be the induced pointed $\underline{H}$-torsor, which can be uniquely lifted to a pointed $\underline{H}$-torsor $(V,y)\to (U,x)$ with the geometrically generic fiber $(V_{\olK},y_{\olK})\to (U_{\olK},x_{\olK})$ . By Lemma \ref{lem:LP^1 diag}, we get a surjective homomorphism of finitely generated $\F_p[H]$-modules
\begin{equation*}
H^1_{\fppf}(V_{\olK},\mu_p)\twoheadrightarrow H^{1}_{\fppf}(V_k,\mu_p)
\end{equation*}
and the $\F_p[H]$-module $M$ can be embedded into $H^1_{\fppf}(V_k,\mu_p)$ as an $\F_p[H]$-submodule. However, as the category $\Mod(\F_p[H])$ is semi-simple, $M$ can be (non-canonically) lifted to an $\F_p[H]$-submodule of $H^1_{\fppf}(V_{\olK},\mu_p)$. By Lemma \ref{lem:count lrq}, this amounts to saying that the pointed $G$-torsor $(P,p)\to (U_k,x_k)$ can be lifted to a pointed $G_{\olK}$-torsor $(P_{\olK},p_{\olK})\to (U_{\olK},x_{\olK})$. This completes the proof. 
\end{proof}

\subsection{The cospecialization map for the tame fundamental group scheme}\label{sec:cosp}

Let $k$ be an algebraically closed field of characteristic $p>0$. Let $R\Def k[[t]]$ and $K\Def\Frac R=k((t))$. Let us fix an algebraic closure $\overline{K}$ of $K$ and $\overline{R}$ the integral closure of $R$ in $\overline{K}$. Let $f\colon X\to S\Def\Spec R$ be a smooth morphism with geometrically connected fibers with an $S$-valued point $x\in X(S)$. We denote by $(X_{\overline{K}},x_{\overline{K}})$ (respectively $(X_k,x_k)$) the geometrically generic fiber (respectively the special fiber) of $(X,x)/S$.

Let $(X^{\tame}_{\overline{K}},x^{\tame}_{\overline{K}})\to (X_{\overline{K}},x_{\overline{K}})$ be the universal pointed $\pi^{\tame}(X_{\overline{K}},x_{\overline{K}})$-torsor. As $\overline{K}$ is algebraically closed, by Proposition \ref{prop:lr alg closed bc}, there exists a profinite linearly reductive $k$-group scheme $\pi^{\tame}(X_{\olK},x_{\olK})_k$ which is unique up to isomorphism such that
\begin{equation*}
\pi^{\tame}(X_{\olK},x_{\olK})\simeq\pi^{\tame}(X_{\olK},x_{\olK})_k\times_k \olK.
\end{equation*}
Then the base change $\pi^{\tame}(X_{\overline{K}},x_{\olK})_{R}\Def\pi^{\tame}(X_{\olK},x_{\olK})_k\times_k R$ is an $R$-model of $\pi^{\tame}(X_{\olK},x_{\olK})$. As $X$ is smooth over $S$, by Proposition \ref{prop:Olsson 1}, the $\pi^{\tame}(X_{\olK},x_{\olK})$-torsor $X^{\tame}_{\overline{K}}\to X_{\overline{K}}$ uniquely extends to a $\pi^{\tame}(X_{\olK},x_{\olK})_{\overline{R}}$-torsor $(X^{\tame}_{\olR},x^{\tame}_{\olR})\to (X_{\olR},x_{\olR})$. By taking the special fiber, we get a pointed $\pi^{\tame}(X_{\olK},x_{\olK})_k$-torsor over $(X_k,x_k)$. As $\pi^{\tame}(X_{\olK},x_{\olK})_k$ is a profinite linearly reductive $k$-group scheme, by the universality of the tame fundamental group scheme $\pi^{\tame}(X_k,x_k)$ for the pointed scheme $(X_k,x_k)$, there exists a unique homomorphism of $k$-group schemes
\begin{equation}\label{eq:cosp tame fgs}
\pi^{\tame}(X_k,x_k)\to\pi^{\tame}(X_{\olK},x_{\olK})_k.
\end{equation}

\begin{definition}\label{def:cosp tame fgs}
We call the map (\ref{eq:cosp tame fgs}) the \textit{cospecialization map} associated with the pointed smooth scheme $(X,x)$ over $S$ and denote it by $\cosp_{(X,x)/S}$ or simply by $\cosp$.
\end{definition}

We also consider a \textit{truncated} version of the cospecialization map. Let $\phi\colon G\to G'$ be a homomorphism of finite linearly reductive $k$-group schemes. As the connected component $G^0$ of the identity is diagonalizable, we have a canonical isomorphism $G^0\xrightarrow{\simeq}\Diag(\X(G^0))$. The $p$-torsion subgroup $\X(G^0)[p]$ of the character group $\X(G^0)$ is stable under the conjugacy action $\pi_0(G)\to\uAut(G^0)\xrightarrow{\simeq}\uAut(\X(G^0))$ by the \'etale quotient $\pi_0(G)$, hence the kernel of the quotient map $G^0\twoheadrightarrow\Diag(\X(G^0)[p])$ is normal in the group scheme $G$. Therefore, we obtain the quotient
\begin{equation*}
G^{\scrD}\Def G/\Ker(G^0\to\Diag(\X(G^0)[p])),
\end{equation*}
which belongs to the class $\scrD$~(cf.\ Definition \ref{def:category D}). 
Similarly for any homomorphism $\phi\colon G\to G'$ of finite linearly reductive group schemes induces a homomorphism $\phi^{\scrD}\colon G^{\scrD}\to {G'}^{\scrD}$. The construction can be generalized to an arbitrary homomorphism between pro-finite linearly reductive group schemes over $k$. Thus, the cospecialization map in Definition \ref{def:cosp tame fgs} induces a homomorphism
\begin{equation}\label{eq:cosp^1 tame fgs}
\cosp_{(X,x)/S}^{\scrD}\colon \pi^{\tame}(X_k,x_k)^{\scrD}\to\pi^{\tame}(X_{\olK},x_{\olK})_k^{\scrD},
\end{equation}   
which we call the \textit{truncated cospecialization map} for the tame fundamental group schemes. 
Now we prove Theorem \ref{thm-int:key}(1).

\begin{prop}\label{prop:cosp}
Let $f\colon X\to S$ be a proper smooth $S$-curve with geometrically connected fibers. Let $D$ be a relatively \'etale Cartier divisor on $X/S$ of degree $n$. We set $U\Def X\setminus \Supp(D)$ and denote by $U_{\olK}$ (respectively by $U_k$) the geometrically generic fiber (respectively the special fiber) of $U/S$. 
Fix an $S$-valued point $x\in U(S)$. Then the truncated cospecialization map (cf.\,(\ref{eq:cosp^1 tame fgs})),
\begin{equation*}
\cosp_{(U,x)/S}^{\scrD}\colon\pi^{\tame}(U_k,x_k)^{\scrD}\to\pi^{\tame}(U_{\olK},x_{\olK})^{\scrD}_k
\end{equation*}
is injective.
\end{prop}

\begin{proof}
Let $\pi^{\tame}(U_k,x_k)^{\scrD}\twoheadrightarrow G$ be an arbitrary finite quotient, which corresponds to some pointed $G$-torsor $(P,p)\to (U_k,x_k)$. By Proposition \ref{prop:LP^1}, the pointed $G$-torsor is lifted to a pointed $G_{\olK}$-torsor over $(U_{\olK},x_{\olK})$. This amounts to saying that the surjective homomorphism $\phi$ factors through $\pi^{\tame}(U_{\olK},x_{\olK})^{\scrD}_k$. As $G$ is arbitrary, this immediately implies that the homomorphism $\cosp^{\scrD}$ is injective.  
\end{proof}

\begin{rem}\label{rem:cosp sm cv}
With the same notation as in Proposition \ref{prop:cosp}, thanks to the base change property for the tame fundamental group schemes~(cf.\ Proposition \ref{prop:bc tame fgs}), the (truncated) cospecialization map induces the map of affine $\olK$-group schemes
\begin{equation*}
\cosp^{\scrD}_{\olK}\colon\pi^{\tame}(U_{k,\olK},x_{k,{\olK}})^{\scrD}\to\pi^{\tame}(U_{\olK},x_{\olK})^{\scrD},
\end{equation*}
where $(U_{k,\olK},x_{k,\olK})$ is the trivial deformation, i.e.\ $U_{k,\olK}=U_{k}\times\olK$ and $x_{k,\olK}=x_k\times\olK$. We also call this map the (truncated) cospecialization map. 
\end{rem}

\begin{thm}\label{thm:cosp}
With the same notation as in Proposition \ref{prop:cosp}, for the pointed $S$-curve $(U,x)$, the following conditions are equivalent to each other.
\begin{enumerate}
\renewcommand{\labelenumi}{(\alph{enumi})}
\item The truncated cospecialization map $\cosp_{(U,x)/S}^{\scrD}$ is isomorphism. 
\item For any $G\in\scrD$~(cf.\ Definition \ref{def:category D}), we have
\begin{equation*}
\#\Hom(\pi^{\tame}(U_k,x_k)^{\scrD},G)=\#\Hom(\pi^{\tame}(U_{\olK},x_{\olK})^{\scrD},G_{\olK}).
\end{equation*}
\item There exists an isomorphism of $k$-group schemes
\begin{equation*}
\pi^{\tame}(U_k,x_k)^{\scrD}\simeq\pi^{\tame}(U_{\olK},x_{\olK})^{\scrD}_k.
\end{equation*}
\end{enumerate}
\end{thm}

\begin{proof}
The implications (a)$\Rightarrow$(c)$\Rightarrow$(b) are clear. It suffices to show the implication (b)$\Rightarrow$(a) holds true. Indeed, by Propositions \ref{prop:finiteness} and \ref{prop:LP^1}, the condition (b) implies that for any $G\in\scrD$, the induced map
\begin{equation*}
\Hom(\pi^{\tame}(U_{\olK},x_{\olK})^{\scrD},G_{\olK})\xleftarrow{\simeq}\Hom(\pi^{\tame}(U_{\olK},x_{\olK})^{\scrD}_k,G)\xrightarrow{\cosp^{\scrD*}}\Hom(\pi^{\tame}(U_k,x_k)^{\scrD},G)
\end{equation*}
is a bijection, which implies that the homomorphism $\cosp^{\scrD}$ itself is an isomorphism. Hence, the condition (a) holds true. This completes the proof.  
\end{proof}

Now we can prove Theorem \ref{thm-int:main1}.

\begin{cor}\label{cor:cosp2}
Let $k_0=\overline{\F}_p$ and $S=\Spec k_0[[t]]=\{s,\eta\}$, where $s$ (respectively $\eta$) is the closed point (respectively generic point) of $S$. Let $X$ be a proper smooth relative $S$-curve of genus $g$ with geometrically connected fibers and $D$  an \'etale relative Cartier divisor on $X$ of degree $n$. We set $U\Def X\setminus \Supp(D)$ and fix an $S$-valued point $x\in U(S)$. Suppose that $U$ is hyperbolic. If $U_{\oeta}$ is not constant~(cf.\ \S\ref{sec:Tamagawa's thm}), then 
$\pi^{\tame}(U_{\oeta},x_{\oeta})$ is not isomorphic to $\pi^{\tame}(U_{s},x_{s})\times_{k_0}k_0(\oeta)$. 
\end{cor}

\begin{proof}
Suppose that $U_{\oeta}$ is a non-constant hyperbolic curve. It suffices to show that $\pi^{\tame}(U_{\oeta},x_{\oeta})^{\scrD}$ is not isomorphic to $\pi^{\tame}(U_{s},x_{s})^{\scrD}\times_{k_0}k_0(\oeta)$. By Theorem \ref{thm:cosp}, we have only to show that the truncated cospecialization map \[
\cosp_{U}^{\scrD}\colon\pi^{\tame}(U_{s},x_{s})^{\scrD}\to\pi^{\tame}(U_{\oeta},x_{\oeta})_{k_0}^{\scrD}
\] 
is not an isomorphism. 
We will adapt the argument in the proof of  \cite[Theorem 8.1]{Ta04}. Suppose that the map $\cosp_{U}^{\scrD}\colon\pi^{\tame}(U_{s},x_{s})^{\scrD}\to\pi^{\tame}(U_{\oeta},x_{\oeta})_{k_0}^{\scrD}$ is an isomorphism.  
As $2-2g-n<0$, there exists a connected pointed $\underline{H}$-torsor $(V_{\oeta},y_{\oeta})\to (U_{\oeta},x_{\oeta})$ where $H$ is a finite group of prime-to-$p$ order such that the normalization $Y_{\oeta}$ of $X_{\eta}$ in $V_{\oeta}$ has genus greater than or equal to $2$. Moreover, we may assume that $Y_{\oeta}\to X_{\oeta}$ is ramified at every point in $\Supp(D_{\oeta})$. Then Grothendieck's specialization theorem for the prime-to-$p$ \'etale fundamental groups
\begin{equation*}
\spe^{{(p')}}_U\colon\pi^{\et}_1(U_{\oeta},x_{\oeta})^{(p')}\xrightarrow{\simeq}\pi_1^{\et}(U_{s},x_{s})^{(p')}
\end{equation*}
(cf.\ (\ref{eq:sp p'})) says that the prime-to-$p$ covering $(V_{\oeta},x_{\oeta})\to (U_{\oeta},x_{\oeta})$ is isomorphic to the geometrically generic fiber of the connected $\underline{H}_{S'}$-torsor $(V,y)\to (U_{S'},x_{S'})$ for some finite extension $S'\to S$. Then by Corollary \ref{cor:tame fgs fet}, we have the commutative diagram of exact sequences 
\begin{equation*}
\begin{xy}
\xymatrix{
1\ar[r]&\pi^{\tame}(V_{s},y_{s})^{\scrD}\ar[r]\ar[d]^{\cosp^{\scrD}_V}&\pi^{\tame}(U_{s},x_{s})^{\scrD}\ar[r]\ar[d]^{\cosp^{\scrD}_{U_{S'}}}_{\simeq}&\underline{H}\ar[r]\ar@{=}[d]&1\\
1\ar[r]&\pi^{\tame}(V_{\oeta},y_{\oeta})_{k_0}^{\scrD}\ar[r]&\pi^{\tame}(U_{\oeta},x_{\oeta})_{k_0}^{\scrD}\ar[r]&\underline{H}\ar[r]&1.
}
\end{xy}
\end{equation*}
This gives an isomorphism $\cosp^{\scrD}_{V}\colon \pi^{\tame}(V_{s},y_{s})^{\scrD}\xrightarrow{\simeq}\pi^{\tame}(V_{\oeta},y_{\oeta})^{\scrD}_{k_0}$ and hence, by Corollary \ref{cor:Olsson 2 open vs proper} (dividing by tame inertia), we get the isomorphism
\begin{equation*}
\cosp^{\scrD}_{Y}\colon \pi^{\tame}(Y_{s},y_{s})^{\scrD}\xrightarrow{\simeq}\pi^{\tame}(Y_{\oeta},y_{\oeta})^{\scrD}_{k_0}.
\end{equation*}
By Propositions \ref{prop:cat C=D} and \ref{prop:fq lr vs et}, this implies that $\pi^{\et}_A(Y_{s})^{\scrC}=\pi^{\et}_A(Y_{\oeta})^{\scrC}$. However, as both the profinite groups $\pi^{\et}_1(Y_{s})^{\scrC}$ and $\pi_1^{\et}(Y_{\oeta})^{\scrC}$ are topologically finitely generated~(cf.\ Corollary \ref{cor:pi^t fin gen}), it turns out that the specialization map $\spe\colon\pi^{\et}_1(Y_{\oeta})^{\scrC}\to\pi_1^{\et}(Y_{s})^{\scrC}$ is an isomorphism~(cf.\ \cite[Lemma 8.4]{Ta04}). By Theorem \ref{thm:Tamagawa C}, we can conclude that $Y_{\oeta}$ is constant. Since $Y_{\oeta}$ is hyperbolic, $\uAut_{k_0(\oeta)}(Y_{\oeta})$ is a constant group scheme. Therefore,  $U_{\oeta}=(Y_{\oeta}/H)\setminus\{\text{ramification locus}\}$ is also constant, which is a contradiction. This completes the proof. 
\end{proof}

\subsection{Reconstruction of numerical invariants}\label{sec:rmk inv}

Let $k$ be an algebraically closed field of characteristic $p>0$. Let $X$ be a proper smooth connected curve over $k$ of genus $g$ and of $p$-rank $\gamma$ (cf.\,Definition \ref{def:p-rank}). Let $S$ be a finite (possibly empty) set of closed points of $X$ and we set $n\Def\#S$. Let $U\Def X\setminus S$. We denote by $\pi^{\tame}(U)$ the isomorphism class of the tame fundamental group scheme $\pi^{\tame}(U,x)$ for some (any) $k$-rational point $x\in U(k)$. We will discuss reconstruction of numerical invariants $(g,n,\gamma)$ from the tame fundamental group scheme $\pi^{\tame}(U)$. This is motivated by the work of Tamagawa \cite{Ta03}, in which the following result is established.

\begin{thm}(cf.\ \cite[Theorem 4.1]{Ta03})\label{thm:Tamagawa03}
For each $i=1,2$, let $p_i$ be a prime number, $k_i$ an algebraically closed field of characteristic $p_i$, $X_i$ a proper smooth connected curve of genus $g_i$ over $k_i$, $S_i$ a finite (possibly empty) set of closed points of $X_i$ with  $n_i\Def\#S_i$ and $U_i\Def X_i\setminus S_i$. If $\pi^{\rmt}_1(U_1)\simeq\pi^{\rmt}_1(U_2)$, then we have $p_1=p_2$, $g_1=g_2$ and $n_1=n_2$ unless $g_i=0$ and $n_i\le 1$ for $i=1,2$. 
\end{thm}

 Note that, as we have $\gamma_i=\Dim_{\F_p}\Hom(\pi^{\et}_1(X_i),\Z/p\Z)=\Dim_{\F_p}\Hom(\pi_1^{\rmt}(U_i),\Z/p\Z)$ 
(cf.\,\S\ref{sec:fq et fg}), the $p$-rank $\gamma_i$ can be easily reconstructed from the tame fundamental group $\pi_1^{\rmt}(U_i)$. The proof involves analysis of the behaviour of the generalized Hasse--Witt invariants for prime-to-$p$ cyclic coverings of curves.

It seems natural to ask whether or not the same result still holds true after replacing the tame fundamental group $\pi^{\rmt}_1(U)$ by the tame fundamental group scheme $\pi^{\tame}(U)$. We assume that $k\Def k_1=k_2$ and set $p\Def p_1=p_2$. We will  freely use the base field $k$ and hence its characteristic $p$. Let $X_i,S_i,U_i,g_i,n_i~(i=1,2)$ be as in Theorem \ref{thm:Tamagawa03}. For each $i=1,2$, let $\gamma_i$ be the $p$-rank of the curve $X_i$ (cf.\,Definition \ref{def:p-rank}).

\begin{lem}\label{lem:rcnst inv}
With the above notation, suppose that we have an isomorphism $\pi^{\tame}(U_1)\simeq\pi^{\tame}(U_2)$ as affine $k$-group schemes. Then we have the following. 
\begin{enumerate}
\renewcommand{\labelenumi}{(\arabic{enumi})}
\item $2g_1+n_1-1+\delta_1=2g_2+n_2-1+\delta_2$.
\item $\gamma_1+n_1-1+\delta_1=\gamma_2+n_2-1+\delta_2$.
\end{enumerate}
Here we set
\begin{equation*}
\delta_i=
\begin{cases}
1&\text{if $n_i=0$},\\
0&\text{if $n_i>0$}.
\end{cases}
\end{equation*}
\end{lem}

\begin{proof}
(1) Let $\ell$ be a prime number with $\ell\neq p$. Then the assertion immediately follows from
\begin{equation*}
2g_i+n_i-1+\delta_i=\Dim_{\F_{\ell}}\Hom(\pi^{\et}_1(U_i)^{(p')},\mu_{\ell}(k))=\Dim_{\F_{\ell}}\Hom(\pi^{\tame}(U_i),\mu_{\ell}).
\end{equation*}

(2) By replacing $\mu_{\ell}$ with $\mu_{p}$ and by using Corollary \ref{cor:Olsson 2 mu_p} or the equation (\ref{eq:p-rank mu_p}), the same argument as in the proof of (1) implies the desired equality. This completes the proof. 
\end{proof}

\begin{ex}\label{ex:rcnst inv}
\begin{enumerate}
\renewcommand{\labelenumi}{(\arabic{enumi})}
\item Let us assume that $n_1=0$, i.e.\ $U_1=X_1$.  Then the condition that $\pi^{\tame}(X_1)\simeq\pi^{\tame}(U_2)$ implies $(g_1,n_1,\gamma_1)=(g_1,0,\gamma_1)=(g_2,n_2,\gamma_2)$ unless $g_i=0$ and $n_i\le 1$ for $i=1,2$. By Lemma \ref{lem:rcnst inv}, the claim holds in the case $n_2=0$. Suppose that $n_2>0$. Let us prove that $g_1=g_2=0$ and $n_2=1$. Indeed, by Corollary \ref{cor:pi^t fin gen}, $\pi^{\et}_1(U_2)^{(p')}$ is a free pro-prime-to-$p$ group of rank $2g_2+n_2-1$. On the other hand, our condition implies that $\pi_1^{\et}(U_1)^{(p')}=\pi^{\tame}(U_1)(k)\simeq\pi^{\tame}(U_2)(k)=\pi_1^{\et}(U_2)^{(p')}$. Hence $\pi^{\et}_1(U_1)^{(p')}$ must be a free pro-prime-to-$p$ group. As $\pi^{\et}_1(U_1)^{(p')}\simeq\Pi^{(p')}_{g_1,0}$~(cf.\ (\ref{eq:Pi g,n})), we must have $g_1\le 1$, in which the tame fundamental group scheme $\pi^{\tame}(X_1)$ is abelian. Hence, so is $\pi^{\tame}(U_2)$. As $n_2>0$, this implies that $g_2=0$. By Lemma \ref{lem:rcnst inv}, we have \begin{equation*}
2g_1=n_2-1=\gamma_1,
\end{equation*}
which is possible only when $g_1=0$ as $\gamma_1\le g_1$. In this case, we have $n_2=1$. 

\item Let us suppose that $g_1=0$, i.e.\ $X_1\simeq\mathbb{P}^1_k$. Then the condition that $\pi^{\tame}(U_1)\simeq \pi^{\tame}(U_2)$ implies that $(g_1,n_1,\gamma_1)=(0,n_1,0)=(g_2,n_2,\gamma_2)$. 
If $n_1=0$ or $n_2=0$, by the first case (1), the assertion is obviously true. Thus, let us assume that $n_1,n_2>0$. Then, by Lemma \ref{lem:rcnst inv}, we have $2g_2+n_2-1=n_1-1=\gamma_2+n_2-1$, hence $2g_2=\gamma_2$. As $\gamma_2\le g_2$, this only happens when $g_2=\gamma_2=0$, in which case we have $n_1=n_2$. 
\end{enumerate}
\end{ex}

Now we prove Theorem \ref{thm-int:main2}.

\begin{thm}\label{thm:rcnst inv}
Suppose that there exists an isomorphism of $k$-group schemes $\pi^{\tame}(U_1)\simeq\pi^{\tame}(U_2)$. Then we have $(g_1,n_1,\gamma_1)=(g_2,n_2,\gamma_2)$ unless $g_i=0$ and $n_i\le 1$ for $i=1,2$.
\end{thm}

\begin{proof}
To ease the notation, we set $\Pi(U_i)\Def\pi^{\tame}(U_i)$ for $=1,2$.  
By Example \ref{ex:rcnst inv}(1),(2), we may assume that $g_i,n_i>0$ for $i=1,2$. In this case, by Lemma \ref{lem:rcnst inv}, we have
\begin{equation*}
2g_1+n_1=2g_2+n_2,\qquad \gamma_1+n_1=\gamma_2+n_2.
\end{equation*}
Hence, it suffices to show that the equality 
$g_1=g_2$ holds true. 
We will apply the results in \cite[\S4]{Ta03}. Fix $i\in \{1,2\}$. For any positive integer $N>0$ with $p\nmid N$, let $U_i(N)\to U_i$ be the finite \'etale Galois covering corresponding to the surjective homomorphism $\Pi(U_i)\twoheadrightarrow \left(\Pi(U_i)^{\et}\right)^{\ab}/N=(\pi_1^{\et}(U_i)^{(p')})^{\ab}/N$. We denote by $X_i(N)$ the normalization of $X_i$ in $U_i(N)$. We set  $n_i(N):=\#\left(X_i(N)\setminus U_i(N)\right)$ for $i=1,2$. Then we have
\begin{equation}\label{eq:tilde gamma vs gamma}
\begin{aligned}
\Dim_{\F_p}H^1_{\fppf}(U_i(N),\mu_p)&=\Dim_{\F_p}H^1_{\fppf}(X_i(N),\mu_p)+n_i(N)-1\\
&=\Dim_{\F_p}H^1_{\et}(X_i(N),\Z/p\Z)+n_i(N)-1,
\end{aligned}
\end{equation}
where the first equality is due to Corollary \ref{cor:Olsson 2 mu_p} and the second one is due to the Serre duality~(\ref{eq:duality X}). Thanks to Corollary \ref{cor:tame fgs fet}, we have 
\begin{equation*}
H^1_{\fppf}(U_i(N),\mu_p)=\Hom\biggl(\Ker\Bigl(\Pi(U_i)\twoheadrightarrow (\Pi(U_i)^{\et})^{\ab}/N\Bigl),\mu_p\biggl).
\end{equation*}
Therefore, the assumption $\Pi(U_1)\simeq\Pi(U_2)$ implies that
\begin{equation}\label{eq:tilde gamma}
\Dim_{\F_p}H^1_{\fppf}(U_1(N),\mu_p)=\Dim_{\F_p}H^1_{\fppf}(U_2(N),\mu_p)
\end{equation}
for any integer $N>0$ with $p\nmid N$. By Tamagawa's theorem~\cite[Corollary 4.11; see also Remark 4.8]{Ta03}, we have
\[
\lim_{N=p^f-1,~f\to\infty}\frac{\Dim_{\F_p}H^1_{\et}(X_i(N),\Z/p\Z)}{N^{2g_i+n_i-1}}=\begin{cases}
g_i-1&\text{if $n_i= 1$},\\
g_i&\text{if $n_i>1$}.
\end{cases}
\]
On the other hand, as
\[
n_i(N)=
\begin{cases}
N^{2g_i+n_i-1}&\text{if $n_i=1$},\\
N^{2g_i+n_i-2}n_i&\text{if $n_i>1$}
\end{cases}
\]
(cf.\,Proof of \cite[Lemma 4.13]{Ta03}), we obtain
\[
\lim_{N=p^f-1,~f\to\infty}\dfrac{n_i(N)-1}{N^{2g_i+n_i-1}}=
\begin{cases}
1&\text{if $n_i=1$},\\
0&\text{if $n_i>1$}.
\end{cases}
\]
Therefore, by (\ref{eq:tilde gamma vs gamma}), we have
\begin{equation}\label{eq:limit av}
\begin{aligned}
&\lim_{N=p^f-1,~f\to\infty}\frac{\Dim_{\F_p}H^1_{\fppf}(U_i(N),\mu_p)}{N^{2g_i+n_i-1}}\\
&=\lim_{N=p^f-1,~f\to\infty}\left(\frac{\Dim_{\F_p}H^1_{\et}(X_i(N),\Z/p\Z)}{N^{2g_i+n_i-1}}+\dfrac{n_i(N)-1}{N^{2g_i+n_i-1}}\right)=g_i.
\end{aligned}
\end{equation}
The equations (\ref{eq:tilde gamma}) and (\ref{eq:limit av}) thus imply that $g_1=g_2$. 
This completes the proof. 
\end{proof}

\begin{rem}
As a consequence of \cite[Theorem 4.1]{Ta03} and \cite[Theorem 8.1]{Ta04}, Tamagawa established the {\it finiteness} theorem \cite[Theorem 8.6]{Ta04}. 
It is natural to ask whether we can deduce our version of the finiteness theorem from Corollary \ref{cor:cosp2} and Theorem \ref{thm:rcnst inv}. 
Thanks to Proposition \ref{prop:fq lr vs et}, this is correct for {\it proper} hyperbolic curves. In the following, we will see that we need more arguments for the general case. 
Let $k_0=\overline{\F}_p$ and $U$ a smooth connected hyperbolic curve over $k_0$. 
Let $\Sigma$ be the set of $k_0$-isomorphism classes of smooth connected curves $U'$ over $k_0$ such that $\pi^{\tame}(U)\simeq\pi^{\tame}(U')$.  
Let $X\Def U^{\rm cpt}$ be the smooth compactification of $U$. Let $g$ be the genus of $X$ and $n\Def \#(X\setminus U)(k_0)$ the cardinality of the complement of $U$. By Corollary \ref{thm:rcnst inv}, we have $\Sigma\subset M_{g,n,k_0}$, where $M_{g,n,k_0}$ is the coarse moduli space over $k_0$ of smooth connected curves of type $(g,n)$. Suppose that $\Sigma$ is an infinite set. Let $C\subset M_{g,n,k_0}$ be an integral $k_0$-curve with $\# (C\cap\Sigma)=\infty$. Let $\overline{\eta}$ be a geometric point over the generic point $\eta$ of $C$. Let $U_{\oeta}$ be the curve corresponding to the geometric point $\oeta\to M_{g,n,k_0}$. Therefore, the problem is reduced to proving that $\pi^{\tame}(U)^{\scrD}\times_{k_0}k_0(\oeta)\simeq\pi^{\tame}(U_{\oeta})^{\scrD}$, which contradicts Corollary \ref{cor:cosp2}. 

In particular, it is necessary that they have the same set of isomorphism classes of finite quotients. 
We denote by $\pi^{\tame}_A(U)^{\scrD}$ (respectively by $\pi^{\tame}_A(U_{\oeta})^{\scrD}$) the set of isomorphism classes of finite quotients of $\pi^{\tame}(U)^{\scrD}\times_{k_0}k_0(\oeta)$ (respectively finite quotients of $\pi^{\tame}(U_{\oeta})^{\scrD}$). Thanks to Proposition \ref{prop:LP^1}, we have $\pi^{\tame}_A(U)^{\scrD}\subseteq\pi^{\tame}_A(U_{\oeta})^{\scrD}$. It suffices to show that $\pi^{\tame}_A(U_{\oeta})^{\scrD}\subseteq\pi^{\tame}_A(U)^{\scrD}$. Let $G\in\pi^{\tame}_A(U_{\oeta})^{\scrD}$ be an arbitrary finite quotient of $\pi^{\tame}(U_{\oeta})^{\scrD}$. By definition, there exists a Nori-reduced $G$-torsor $P_{\oeta}\to U_{\oeta}$ (cf.\,\ref{para:Nori-reduced}).  
The problem is that the Nori-reducedness does not mean that $P_{\oeta}$ is reduced in general. If $P_{\oeta}$ is reduced, we can conclude that $G\in\pi_A^{\tame}(U)^{\scrD}$ as follows.  
After shrinking $C$ if necessary, there exist a finite flat morphism between integral $k_0$-curves $C'\to C$ and a $G$-torsor $P\to U_{C'}$ over a smooth model  $U_{C'}$ of $U_{\oeta}$ over $C'$ such that $P|_{U_{\oeta}}\simeq P_{\oeta}$ as $G$-torsors over $U_{\oeta}$. 
If $P_{\oeta}$ is reduced, according to \cite[Lemma 0578]{stack}, by shrinking $C$ if necessary, we may assume that  $P\to C'$ has geometrically reduced fibers. If $s\in C'(k_0)$ is a closed point of $C'$ whose image in $C$ belongs to $\Sigma$, we obtain a $G$-torsor $P|_{U_{s}}\to U_s$ with $P|_{U_s}$ reduced and connected. This implies that $G\in\pi_A^{\tame}(U_s)^{\scrD}=\pi_A^{\tame}(U)^{\scrD}$. 
\end{rem}


\end{document}